\documentclass[11pt]{sat}

\usepackage[latin1]{inputenc}
\usepackage{latexsym}
\usepackage{amsmath}

\title{M\"{u}ntz Type Theorems I}
\def\shorttitle{M\"{u}ntz Type Theorems}

\author{J. M. Almira}
\def\shortauthor{J. M. Almira}

\def\versiondate{12 October, 2007}

\def\keywords{} 

\newtheorem{theorem}{Theorem}

\newtheorem{corollary}[theorem]{Corollary}

\newtheorem{lemma}[theorem]{Lemma}

\newtheorem{proposition}[theorem]{Proposition}
\newtheorem{remark}[theorem]{Remark}

\def\Proof#1. {\par
                      \ifdim\lastskip<15pt
                      \removelastskip\penalty-200
                      \vskip15pt plus3pt minus3pt
                      \fi
                       {\def\a{#1}
                       \ifx\a\empty
                       {\noindent\bf Proof.}
                       \else
                       {\noindent\bf Proof of #1.}
                       \fi}\enspace}
\def\restr#1{\,\vrule\,\lower1ex\hbox{$#1$}}

\def\dword#1{{\bf#1}}
\font\Blbb=msbm10 
\def\mathbb#1{\hbox{{\Blbb#1}}}
\let\Bbb\mathbb
\def\dd{\,{\rm d}}  
\def\ee{{\rm e}}  
\def\ii{{\rm i}}  
\def\endproof{\hfill$\Box$\bigskip}
\def\Tchebychev{Chebyshev}
\def\floor#1{\lfloor#1\rfloor} 

\input colordvi

\def\startpagenumber{152}
\def\volumenumber{3}
\def\year{2007}

\newcommand{\beginddoc}{
\begin{document}
\maketitle
\begin{abstract}
In this paper, we concentrate our attention on the M\"{u}ntz
problem in the univariate setting and for the uniform norm.

\vskip1pt MSC: 41-01, 41-02 \ifx\keywords\empty\else\vskip1pt
keywords: \keywords\fi
\end{abstract}
\insert\footins{\scriptsize
\medskip
\baselineskip 8pt \leftline{Surveys in Approximation Theory}
\leftline{Volume \volumenumber, \year.
pp.~\thepage--\pageref{endpage}.} \leftline{Copyright \copyright\
2006 Surveys in Approximation Theory.} \leftline{ISSN 1555-578X}
\leftline{All rights of reproduction in any form reserved.}
\smallskip
\par\allowbreak}
\tableofcontents}
\renewcommand\rightmark{\ifodd\thepage{\it \hfill\shorttitle\hfill}\else {\it \hfill\shortauthor\hfill}\fi}
\markboth{{\it \shortauthor}}{{\it \shorttitle}} \markright{{\it
\shorttitle}}
\def\endddoc{\label{endpage}\end{document}}
\date{{\small \versiondate}}
\setlength\oddsidemargin{0pc} \setlength\evensidemargin{0pc}
\setlength\topmargin{0in} \setlength\textwidth{6.5in}
\setlength\textheight{8.6in}
\beginddoc



\section{Introduction}

In his seminal paper \cite{bernstein1912}
 of 1912, the Russian mathematician S. N.
Bernstein (one of the greatest approximation theorists of the last
century) asked under which conditions on an increasing sequence
$\Lambda =(0=\lambda _{0}<\lambda _{1}<\cdots)$ one can guarantee
that the vector space
\begin{equation}\label{densidadmuntz}
\Pi(\Lambda) := {\rm span}\{x^{\lambda_k}: k=0,1,\ldots\}
\end{equation}
spanned by the monomials $x^{\lambda _{k}}$ is \ a dense subset of
$C[0,1]$. He specifically proved that the condition
\[
\sum_{\lambda _{k}>0}\frac{1+\log \lambda _{k}}{\lambda
_{k}}=\infty
\]
is necessary and the condition
\[
\lim_{k\to \infty} \frac{\lambda_{k}}{k \log k} = 0
\]
is sufficient, and conjectured that a necessary and sufficient
condition to have $\overline{\Pi(\Lambda)}= C[0,1]$ is
\[
\sum_{k=1}^{\infty }\frac{1}{\lambda _{k}}=\infty.
\]

This conjecture was proved by M\"{u}ntz \cite{muntz} in 1914. In
his proof, he
used the method of Gram determinants to compute the distance of $x^{\lambda }$%
from $\Pi(\lambda_0,\ldots,\lambda_n))_2$ in 
the $L^{2}(0,1)$--metric. The determinants that appear in this
problem are of the form
\[\det (1/(1+a_{i}+a_{j}))_{0\leq
i,j\leq n},\] and their explicit expression was obtained in the
19th century by Cauchy.

For the sake of clarity, let us give a precise formulation of the
classical M\"{u}ntz Theorem.

\begin{theorem}[M\"{u}ntz, 1914]\label{muntzproof} Let\textbf{\ }$\Lambda =(\lambda _{i})_{k=0}^{\infty
}$, $0=\lambda _{0}<\lambda _{1}<\cdots $, be an increasing
sequence of non-negative real numbers. Then
$\Pi (\Lambda )={\rm span}\{x^{\lambda _{k}}:k=0,1,\ldots\}$, the
\dword{M\"{u}ntz space associated to $\Lambda$}, is a dense subset
of $C[0,1]$ if and only if
\begin{equation}\label{harmonic}
\sum_{k=1}^{\infty }\frac{1}{\lambda _{k}}=\infty .
\end{equation}
\end{theorem}

This is a beautiful theorem because it connects a topological
result (the density of a certain subset of a functional space)
with an arithmetical one (the divergence of a certain harmonic
series). Many people might well have been drawn to this result
because of its beauty. Another reason to be interested in
M\"{u}ntz' theorem is that the original result not only solves a
nice problem but also opens the door to many new interesting
questions. For example, one is tempted to change the space of
continuous functions $C[0,1]$ to other function spaces such as
$L^p(a,b)$, or to consider the analogous problem in several
variables, on complex domains, on intervals away from the origin,
for more general exponent sequences, for polynomials with integral
coefficients, etc. As a consequence, many proofs (and
generalizations) of the theorem have been produced.

In this paper, we concentrate our attention on the M\"{u}ntz
problem in the univariate setting and for the uniform norm.
Moreover, we do not include any results about the rate of
convergence to zero of the errors of best (uniform) approximation
using M\"{u}ntz polynomials. On the other hand, we do provide
proofs in great detail, and we promise to write a second paper
where we plan to treat several advanced topics, including the
M\"{u}ntz Theorem for complex domains, M\"{u}ntz-Jackson theorems,
M\"{u}ntz type theorems for approximation with polynomials with
integral coefficients, the $p$-adic M\"{u}ntz theorems, and the
M\"{u}ntz Theorem for rational functions.

Let us return to a discussion of this paper. We devote 
Section 2 to the classical M\"{u}ntz Theorem. In particular, we give
several proofs of this result, showing how the  M\"{u}ntz problem
is connected to many apparently different branches of mathematics.
In Section 3, we focus our attention on the so called
Full M\"{u}ntz Theorem, i.e., we study the density of ${\rm
span}\{x^{\lambda_k}\}_{k=0}^{\infty}$ in $C(K)$, where $K$
denotes a compact subset of $[0,\infty)$, for arbitrary sequences
of exponents
$(\lambda_k)_{k=0}^{\infty}$ and characterize the (uniform)
closure of ${\rm span}\{x^{\lambda_k}\}_{k=0}^{\infty}$ in the
nondense case.

\section{The classical M\"{u}ntz theorem}
\subsection{M\"{u}ntz theorem: the original proof with a modification by O. Sz\'{a}sz}
The original proof by M\"{u}ntz of Theorem \ref{muntzproof}, and
that which remains essentially the standard proof that
you may find in many introductory textbooks on approximation theory, is based on an estimation of the errors $%
E(x^{q},\Pi (\Lambda _{n}))$, where $\Lambda _{n}:=(\lambda
_{k})_{k=0}^{n}$
 and
\[
E(x^{q},\Pi (\Lambda _{n})):=\inf_{p\in \Pi (\Lambda
_{n})}\|x^{q}-p(x)\|_{[0,1]}
\]
is the error of best approximation, with respect to the uniform
norm in $[0,1]$, to $x^{q}$ when we take as the set of
approximants the space $\Pi (\Lambda _{n})$. It is clear that $\Pi
(\Lambda )$ is dense in $C[0,1]$ if and only if  for all $q\in
\mathbb{N}$, $E(x^{q},\Pi (\Lambda _{n}))$ converges to zero as
$n$ tends to infinity (this is a consequence of the Weierstrass
Approximation Theorem). So, how do we produce a reasonable
estimate for $E(x^{q},\Pi (\Lambda _{n}))$?

If we use the $L_{2}(0,1)$-norm, then we can explicitly compute
the errors
\[
E(x^{q},\Pi (\Lambda _{n}))_{2}=\inf_{p\in \Pi (\Lambda
_{n})}\left\| x^{q}-p(x)\right\| _{2}
\]
since $L_{2}(0,1)$ is a Hilbert space. In fact, if we denote by $%
G(f_{1},\ldots,f_{n})$ the Gram determinant associated with a
linearly independent sequence $(f_{1},\ldots,f_{n})$ of elements
in a Hilbert space $H$ with inner product $(\cdot ,\cdot )$,
\[
G(f_{1},\ldots,f_{n})=\det \left(
\begin{array}{clc}
(f_{1},f_{1}) & \cdots & (f_{1},f_{n}) \\
\vdots & \ddots & \vdots \\
(f_{n},f_{1}) & \cdots & (f_{n},f_{n})
\end{array}
\right),
\]
then it is well known \cite[Theor. 8.7.4.]{davis} that
\[
E(g,V_{n})_{H}=\inf_{v\in V_{n}}\left\| g-v\right\| _{H}=\sqrt{\frac{%
G(g,f_{1},\ldots,f_{n})}{G(f_{1},\ldots,f_{n})}}
\]
holds for all $g\notin V_{n}={\rm span}\{f_{1},\ldots,f_{n}\}$.
From this follows (for $\lambda _{0}>-1/2$) the formula
\[
E(x^{q},\Pi (\Lambda _{n}))_{2}=\frac{1}{\sqrt{2q+1}}\prod_{k=0}^{n}\frac{%
\left| q-\lambda _{k}\right| }{q+\lambda _{k}+1}
\]
for all $q>-1/2$, since

\[
G(x^{\lambda_{1}},\ldots,x^{\lambda_{n}})=\det \left(
\frac{1}{\lambda_{i}+\lambda_{j}+1} \right)_{0\leq i,j\leq
n}=\frac{\prod_{i>j}^n(\lambda_i-\lambda_j)^2}{\prod_{i,j=1}^n(\lambda_i+\lambda_j+1)}.
\]
The Cauchy determinant $\det \left(
\frac{1}{\lambda_{i}+\lambda_{j}+1} \right)_{0\leq i,j\leq n}$ is
a particular case of $\det \left( \frac{1}{a_{i}+b_{j}}
\right)_{1\leq i,j\leq n}$ and the argument to compute a closed
expression for this determinant is quite similar to the classical
argument to compute a Vandermonde determinant. Thus, it consists
in considering both sides of the identity
\begin{equation}\label{cauchy}
\prod_{i=1}^n\prod_{j=1}^n(a_i+b_j)\det\left(
\frac{1}{a_{i}+b_{j}} \right)_{1\leq i,j\leq
n}=\prod_{i=1}^n\prod_{j=1}^{i-1}(a_i-a_j)(b_i-b_j)
\end{equation}
as polynomials in the variables $a_i,b_j$, each time taking into
account just one of these variables, and using the zero properties
of algebraic polynomials of one variable to prove that both
expressions are the same. A detailed proof of $(\ref{cauchy})$ can
be found in \cite[page. 74]{gurariy} or \cite[page 268]{davis}.

From here it is not difficult to prove that
\[
\lim_{n\rightarrow \infty }E(x^{q},\Pi (\Lambda _{n}))_{2}=0\; \text{ for all }%
q\in \mathbb{N}\quad\text{ if and only if }\quad \sum_{k=1}^{\infty }\frac{1}{\lambda _{k}}%
=\infty.
\]
Hence $\Pi (\Lambda )$ is dense in $L_{2}(0,1)$ if and only if $%
\sum_{k=1}^{\infty }1/\lambda _{k}=\infty $. This clearly implies
the necessity of the condition $\sum_{k=1}^{\infty }1/\lambda _{k}%
=\infty $ \ to guarantee the density of $\Pi (\Lambda )$ in
$C[0,1]$.

In order to guarantee the sufficiency of condition
(\ref{harmonic}) to the claim that
$\overline{\Pi(\Lambda)}=C[0,1]$, M\"{u}ntz used Fej\'{e}r's
theorem on summation of Fourier series, but his proof is too
complicated to be reproduced here. In 1916 Otto Sz\'{a}sz extended
M\"{u}ntz's theorem in the sense that he was able to prove the
result also for certain special sequences of
complex numbers $(\lambda _{k})_{k=0}^{\infty }$ as exponents (see
\cite{szasz1}). Furthermore, he simplified the final step of
M\"{u}ntz's proof, showing that the result in $L^{2}(0,1)$ implies
the same result in $C[0,1]$. This follows from the inequality
\begin{eqnarray*}
\left| x^{q}-\sum_{k=1}^{n}a_{k}x^{\lambda _{k}}\right| &=&\left|
\int_{0}^{x}\left( qt^{q-1}-\sum_{k=1}^{n}a_{k}\lambda
_{k}t^{\lambda_{k}-1}\right) \dd t\right| \\[10pt]
 &\leq
&\int_{0}^{1}\left| qt^{q-1}-\sum_{k=1}^{n}a_{k}\lambda_{k}t^{\lambda_{k}-1}\right| \dd t \\[10pt]
&\leq &\left[ \int_{0}^{1}\left|
qt^{q-1}-\sum_{k=1}^{n}a_{k}\lambda_{k}t^{\lambda _{k}-1}\right| ^{2}\dd t\right] ^{1/2} \\[10pt]
&=&\left\| qx^{q-1}-\sum_{k=1}^{n}a_{k}\lambda
_{k}x^{\lambda_{k}-1}\right\| _{L^{2}[0,1]}
\end{eqnarray*}
which holds for all $x\in [0,1]$. In other words,
\begin{equation}\label{trucoszasz}
\left\|
x^{q}-\sum_{k=1}^{n}a_{k}x^{\lambda_{k}}\right\|_{\mathbf{C}[0,1]}
\leq
\left\|qx^{q-1}-\sum_{k=1}^{n}a_{k}\lambda_{k}x^{\lambda_{k}-1}\right\|
_{L^{2}[0,1]}.
\end{equation}
\endproof

For historical reasons, we include here (without proof) the
precise statement of Sz\'{a}sz's theorem.

\begin{theorem}[Sz\'{a}sz, 1916]\label{szaszteorema} Let
\[ C([0,1],\mathbb{C}):=\{f:[0,1]\to\mathbb{C}: f\text{ continuous}\}
\]
be the space of continuous complex-valued functions defined on
$[0,1]$ and assume that the M\"{u}ntz polynomials have complex
coefficients and complex exponents, so that for
$\Lambda=(\lambda_k)_{k=0}^{\infty}\subset \mathbb{C}$ we set
$\Pi_{\mathbb{C}}(\Lambda):={\rm
span}_{\mathbb{C}}\{x^{\lambda_k}\}_{k=0}^{\infty}$. If
$\lambda_0=0$ and ${\rm Re}(\lambda_k)>0$ for all $k>0$, then
$\Pi_{\mathbb{C}}(\Lambda)$ is a dense subset of ${\rm
C}([0,1],\mathbb{C})$ whenever
\begin{equation}\label{condicionszasz}
\sum_{k=1}^{\infty}\frac{{\rm
Re}(\lambda_k)}{1+|\lambda_k|^2}=\infty.
\end{equation}
Moreover, if
\[\sum_{k=1}^{\infty}\frac{1+{\rm Re}(\lambda_k)}{1+|\lambda_k|^2}<\infty,\]
then $\Pi_{\mathbb{C}}(\Lambda)$ is not a dense subset of ${\rm
C}([0,1],\mathbb{C})$. In particular, if
\[\liminf_{k\to\infty}{\rm Re}(\lambda_k)>0\] then $\Pi_{\mathbb{C}}(\Lambda)$
 is
a dense subset of ${\rm C}([0,1],\mathbb{C})$ if and only if
$(\ref{condicionszasz})$ holds.
\end{theorem}

Note that Sz\'{a}sz's theorem is not conclusive for all cases. For
example, the sequence $\lambda_k=\frac{1}{k}+ \ii\sqrt{k}$
satisfies
\[
\sum_{k=1}^{\infty}\frac{{\rm
Re}(\lambda_k)}{1+|\lambda_k|^2}<\infty\qquad \text{ and } \qquad
\sum_{k=1}^{\infty}\frac{1+{\rm
Re}(\lambda_k)}{1+|\lambda_k|^2}=\infty.\]


\subsection{M. von Golitschek's constructive proof of the M\"{u}ntz theorem}
In this subsection, we present another proof of the fact that
condition (\ref{harmonic}) is sufficient if the relation $\lambda
_{k}\to\infty $ holds. To do this, we follow a very nice proof
published by M. von Golitschek in \cite {golitschekproof}, which
has two distinct advantages with respect to the majority of known
proofs of the same result. It is both constructive and short.

The idea is to define, for each $q>0$, a concrete sequence of
approximants to $x^{q}$, $(P_{n})_{n=0}^{\infty }\subset \Pi
(\Lambda )$ and to prove that $Q_n(x)=x^q-P_n(x)$ converges to
zero uniformly on $[0,1]$. So, we set $Q_{0}(x):=x^{q}$, and, for
$n=1,2,\ldots$, if we already know that
$$Q_{n-1}(x)=x^q-\sum_{k=1}^{n-1}a_{k,n-1}x^{\lambda_k}$$ with some
coefficients $a_{k,n-1}$, then let
\begin{eqnarray*}
Q_{n}(x) &:=&(\lambda _{n}-q)x^{\lambda
_{n}}\int_{x}^{1}Q_{n-1}(t)t^{-(1+\lambda _{n})}\dd t \\[10pt]
&=&(\lambda _{n}-q)x^{\lambda _{n}}\int_{x}^{1}\left(
t^{q}-\sum_{k=1}^{n-1}a_{k,n-1}t^{\lambda _{k}}\right)
t^{-(1+\lambda
_{n})}\dd t \\[10pt]
&=&(\lambda _{n}-q)x^{\lambda _{n}}\int_{x}^{1}\left(
t^{q-(1+\lambda _{n})}-\sum_{k=1}^{n-1}a_{k,n-1}t^{\lambda
_{k}-(1+\lambda _{n})}\right) \dd t
\\[10pt]
&=&(\lambda _{n}-q)x^{\lambda _{n}}\left[ \frac{t^{q-\lambda _{n}}}{%
q-\lambda _{n}}-\sum_{k=1}^{n-1}a_{k,n-1}\frac{t^{\lambda _{k}-\lambda _{n}}%
}{(\lambda _{k}-\lambda _{n})}\right] _{x}^{1} \\[10pt]
&=:&x^{q}-\sum_{k=1}^{n}a_{k,n}x^{\lambda _{k}},
\end{eqnarray*}
hence $P_{n}(x)= \sum_{k=1}^{n}a_{k,n}x^{\lambda _{k}}$.

We only need to prove that $\left\| Q_{n}\right\| _{C[0,1]}$
converges to zero as $n\rightarrow \infty $. Now $\left\|
Q_{0}\right\| _{C[0,1]}=1$, and for all $n\in \mathbb{N}$ we get
from the inequality
\[\lambda x^\lambda(1-x)<1\qquad\mbox{for all $x\in (0,1)$ and $\lambda>0$}\]
that
\[
\left\| Q_{n}\right\| _{C[0,1]}\leq \left| 1-\frac{q}{\lambda _{n}}%
\right| \left\| Q_{n-1}\right\| _{C[0,1]}.
\]
Hence
\[
\left\| Q_{n}\right\| _{C[0,1]}\leq \prod_{k=0}^{n}\left| 1-\frac{q%
}{\lambda _{n}}\right| \rightarrow 0,\qquad n\rightarrow \infty .
\]
\endproof

\subsection{Measure-theoretic focus}
The classical M\"{u}ntz Theorem can be formulated in terms of
measures. We explain here the way in which this formulation is
attained and we give a proof, due to W. Feller \cite{feller}, of
the `only if' part of the result based on measure theoretical
considerations. For pedagogical reasons, we postpone Feller's
proof of the `if' part to the next subsection.

Let us assume that $\Lambda=(\lambda_k)_{k=1}^{\infty}$ is an
increasing sequence of positive real numbers and, to avoid
problems with the origin, let us also assume in this subsection that
the functions we want to approximate vanish at the origin. In this
case, we can rephrase the classical M\"{u}ntz Theorem as follows:
{\sl The space $\Pi(\Lambda)$ is a dense subset of
$C_0[0,1]:=C[0,1]\cap \{f: f(0)=0\}$ if and only if
 $\sum_{k=1}^{\infty}1/\lambda_k=\infty$}. When dealing
 with the problem of the density of
certain linear subspaces of  a Banach space, it is a quite natural
to use the Hahn-Banach Theorem in the following way: {\sl if $Y$ is a
closed subspace of the Banach space $X$, then $Y\neq X$ if and
only if there exists a bounded linear functional $L\in X^*$ such
that $L\neq 0$ and $L|_Y=0$}. Thus, taking $X=C_0[0,1]$ and
$Y=\overline{\Pi(\Lambda)}$, we have that
$\overline{\Pi(\Lambda)}\neq C_0[0,1]$ if and only if there exists
an  $L\in C^*_0[0,1]\setminus\{0\}$ satisfying
$L(x^{\lambda_k})=0$ for all $k=1,2,\ldots$. The dual of
$C_0[0,1]$ is characterized (by the Riesz Representation Theorem)
as follows: $L\in C^*_0[0,1]$ if and only if
\begin{equation}
L(f)=\int_0^1f(t)\dd \mu(t)
\end{equation}
for a certain finite signed Borel measure $\mu$ on $(0,1]$.
Moreover, we know (by the Weierstrass Approximation Theorem) that
algebraic polynomials that vanish at $0$ form a dense subspace of
$C_0[0,1]$. Hence a new formulation of the classical M\"{u}ntz
Theorem is given as follows.

\begin{theorem}[Classical M\"{u}ntz Theorem in terms of
Measures]\label{muntzmeasures} Let us assume that
$(\lambda_k)_{k=1}^{\infty}$ is an increasing sequence of positive
real numbers and let us define, for each finite signed Borel
measure $\mu$ supported on $(0,1]$, the function
\begin{equation}\label{fmu}
f(z):=\int_0^1 t^z \dd \mu(t).
\end{equation}
Then the following claims are equivalent:
\begin{itemize}
\item[$(a)$] $ \sum_{k=1}^{\infty}\frac{1}{\lambda_k} <\infty $

\item[$(b)$] There exists a $($non-zero$)$ finite signed Borel
measure $\mu$ on $(0,1]$ such that $f(\lambda_k)=0$ for all $k\geq
1$ (where $f$ is given by $(\ref{fmu})$).
\end{itemize}
\end{theorem}

Let us prove that $(a)\Rightarrow (b)$. Given the measure $\mu$ we
make the change of variable $t=\ee^{-s}$ which transforms the
interval $(0,1]$ onto the interval $[0,\infty)$, the measure $\mu$
into another measure $m$ on $[0,\infty)$ and the expression
(\ref{fmu}) to the new formula:
\begin{equation}\label{fmudos}
f(z)=\int_0^\infty \ee^{-zs} \dd m(s).
\end{equation}
Hence we should prove that under condition $(a)$ there exists a
finite signed  Borel measure $m$ supported on $[0,\infty)$ such
that the function $f$ given by (\ref{fmudos}) is not identically
zero but satisfies $f(\lambda_k)=0$ for all $k\geq 1$. We present
a proof whose order is reversed: we first define a function $f$
that satisfies $f(\lambda_k)=0$ for all $k\geq 1$ and then we
prove that this function admits an expression of the form
(\ref{fmudos}).

Set, with $\eta>0$,
\[
f(t):=\frac{1}{(1+\eta+t)^2}\prod_{k=1}^{\infty}\frac{\lambda_k-t}{\lambda_k+2\eta+t}\;.
\]
Obviously, $(a)$ guarantees the convergence of the infinite
product defining $f$. This convergence is uniform and absolute on
compact subsets of
$\mathbb{C}\setminus\{-\lambda_k-2\eta\}_{k=1}^{\infty}$. In
particular, $f$ is well defined on $[0,\infty)$ and vanishes on
the sequence $(\lambda_k)_{k=1}^{\infty}$.

Let us define
\[
f_0(t):=\frac{1}{(1+\eta+t)^2}\qquad \text{ and }\qquad
f_k(t):=\frac{\lambda_k-t}{\lambda_k+2\eta+t}f_{k-1}(t),\text{ for
} k=1,2,\ldots.
\]
It is clear that
\begin{equation*}
f_0(t)=\int_0^\infty s\ee^{-(1+t+\eta)s}\dd s=\int_0^\infty
\ee^{-ts}u_0(s)\dd s,
\end{equation*}
where $u_0(s):=s\ee^{-(1+\eta)s}$. Let us assume that, for all
$k<n$, the function $f_k$ admits an expression of the form
\[
f_k(t)=\int_0^\infty \ee^{-ts}u_k(s)\dd s\qquad
\]
with $u_k(0)=0$ (which we know to be true for $k=0$). We will prove that this is
then also the case for $k=n$. In fact, taking into account the
recursive definition of $f_n$,
\[
f_n(t)=\frac{\lambda_n-t}{\lambda_n+2\eta+t}f_{n-1}(t)=\frac{\lambda_n-t}{\lambda_n+2\eta+t}
\int_0^\infty \ee^{-ts}u_{n-1}(s)\dd s,
\]
and, integrating by parts, we note that
\[
tf_{n-1}(t)=\int_0^\infty \ee^{-ts}u_{n-1}'(s)\dd s,
\]
and
\[
f_n(t)=\int_0^\infty \ee^{-ts}u_n(s)\dd s,
\]
where $u_n$ is the solution of the initial value problem

\begin{equation*}
\left\{%
\begin{array}{lll}
 u_n'+u_{n-1}' & = &  \lambda_n u_{n-1}-(\lambda_n+2\eta)u_n\\[10pt]
  u_n(0) & = & 0. \\
\end{array}
\right.
\end{equation*}
Moreover, it is possible to check that under these conditions the
solution $u_n$ satisfies $\lim_{t\to \infty}u_n(t)=0$. Let us now
multiply both sides of the differential equation defining $u_n$ by
$(u_n+u_{n-1})$. We get
\begin{eqnarray*}
\frac{1}{2}[(u_n+u_{n-1})^2]' &=& (u_n'+u_{n-1}')(u_n+u_{n-1})
 = (u_n+u_{n-1})( \lambda_n u_{n-1}-(\lambda_n+2\eta)u_n)\\
 &=& \lambda_n(u_{n-1}^2-u_n^2)-\eta(2 u_n^2+2u_nu_{n-1})\\[5pt]
 &\leq& (\lambda_n+\eta)(u_{n-1}^2-u_n^2).
\end{eqnarray*}
Hence, for all $h\in (0,\infty)$,
\[
\frac{1}{2}(u_n(h)+u_{n-1}(h))^2=\int_{0}^h
\frac{1}{2}[(u_n+u_{n-1})^2]'\dd t\leq
(\lambda_n+\eta)\int_{0}^h(u_{n-1}^2(t)-u_n^2(t))\dd t,
\]
therefore
\[
\int_0^{\infty}u_n^2(t)\dd t\leq \int_0^{\infty}u_{n-1}^2(t)\dd
t,\qquad n=1,2,3,\ldots.
\]
Taking into consideration the convergence of $f_n$ to $f$ and the
weak sequential compactness of the unit ball of $L^2(0,\infty)$,
we conclude that there exists a function $u\in L^2(0,\infty)$ such
that
\[
f(t)=\int_0^{\infty}\ee^{-ts}u(s)\dd s \qquad \text{ for all }
\quad t\geq 0.
\]
Moreover, the same arguments we have used for $f$ should work with
$f^*(t):=f(t-\eta)$ (choose $\lambda_k^*=\lambda_k+\eta$ and
$\eta^*=0$ instead of the old values $\lambda_k$ and $\eta$).
Hence
\[
f^*(t)=\int_0^{\infty}\ee^{-ts}u^*(s)\dd s \qquad \text{ for all }
\quad t\geq 0,
\]
with $u^*(s):=\ee^{\eta s}u(s)\in L^2(0,\infty)$. Of course, this
implies that
\[\int_0^{\infty}|u(t)|\dd t<\infty,\]
 so that $f$ is of
the form (\ref{fmudos}) with $m$ the signed measure that has
density $u$.

\subsection{Two proofs of the `if' part based on the use of
divided differences} One of the first things we observe when
studying the M\"{u}ntz Theorem is that the necessity of condition
(\ref{harmonic}) and its sufficiency are two facts of a quite
different nature, so that the proofs of the `if' part and the
`only if' part of the M\"{u}ntz Theorem are usually independent.
Hence it is tempting to present new proofs for each one of these
parts in terms of one's own interest in the subject.

In this subsection, we will explain two proofs of the `if' part of
the classical M\"{u}ntz Theorem, both based on the use of divided
differences.


The first proof is due to W. Feller \cite{feller}. It is a natural
continuation of the proof given in the previous subsection. It uses
the strong connection between divided differences and completely
monotone functions and certain results from functional analysis
and measure theory. The second proof, by Hirschman and Widder
\cite{widder} and Gelfond \cite{gelfond}, uses divided differences
to construct an adequate generalization of Bernstein polynomials
with the property that the new polynomials only depend on the
powers $\{x^{\lambda_k}\}_{k=1}^{\infty}$.

First, we would like to recall the definition of divided
differences. Given a function $f$ and a subset
$\{x_k\}_{k=0}^{\infty}$ of its domain, we define the divided
differences of $f$ with respect to the nodes
$\{x_k\}_{k=0}^{\infty}$ recursively:
\[
f[x_k]:=f(x_k),\text{ and }
f[x_{i_0},\ldots,x_{i_n}]:=\frac{f[x_{i_0},x_{i_1},\ldots,x_{i_{n-1}}]-f[x_{i_1},x_{i_2},\ldots,x_{i_{n}}]}{x_{i_0}-x_{i_n}}.
\]
These numbers can be characterized in many ways. One of their main
properties is that they are the coefficients in the Newton
representation of the Lagrange interpolation polynomial of $f$ at
the nodes $\{x_0,x_1,\ldots,x_n\}$. More precisely, if
$P_n(x)=a_0+a_1x+\cdots+a_nx^n$ is the unique polynomial of degree
$\leq n$ that satisfies $P(x_k)=f(x_k)$, $k=0,1,\ldots n$, then
\begin{equation}\label{newtonpolynomial}
P_n(x)=f[x_0]+f[x_0,x_1](x-x_0)+\cdots+f[x_0,x_1,\ldots,x_n](x-x_0)(x-x_1)\cdots(x-x_{n-1}),
\end{equation}
and this characterizes the values
\[
f[x_0], f[x_0,x_1],\ldots, f[x_0,\ldots,x_n].
\]
An easy consequence of (\ref{newtonpolynomial}) is that, for all
$x$, the error $R_n(x):=f(x)-P_n(x)$ is given by
\begin{equation}\label{errornewton}
R_n(x)=f[x,x_0,x_1,\ldots,x_n](x-x_0)(x-x_1)\cdots(x-x_n).
\end{equation}
Moreover, for functions $f$ in $C^{(n)}(I)$, where
$I=[\min\{x_i\}_{i=0}^n,\max\{x_i\}_{i=0}^n]$, taking into account
that $R_n(x_i)=0$ for $i=0,1,\cdots,n$, we conclude that $R_n'(x)$
has at least $n$ zeros in the interval $I$, $R''_n(x)$ has at
least $n-1$ zeros therein, etc., so that $R_n^{(n)}(\tau)=0$ for a
certain value $\tau\in I$. Now,
\[
R_n^{(n)}(\tau)=f^{(n)}(\tau)-n!f[x_0,x_1,\cdots,x_n],
\]
so that for a function $f$ that is sufficiently many times
differentiable, the divided differences satisfy
\begin{equation}\label{diferenciasdivididasderivadas}
f[x_0,x_1,\ldots,x_n]=\frac{1}{n!}f^{(n)}(\tau)
\end{equation}
for a certain $\tau\in
I=[\min\{x_i\}_{i=0}^n,\max\{x_i\}_{i=0}^n]$.

Finally, it is also useful to note that:
\begin{equation}\label{diferenciasexpandidas}
f[x_0,x_1,\ldots,x_n]=\sum_{k=0}^n\frac{f(x_k)}{(x_k-x_0)(x_k-x_1)\cdots(x_k-x_{k-1})(x_k-x_{k+1})\cdots(x_k-x_n)}.
\end{equation}
This follows from  the fact (see (\ref{newtonpolynomial}))  that $f[x_0,\ldots,x_n]$ is
the coefficient of $x^n$ in the power form of the interpolating
polynomial, and the Lagrange expression of this polynomial.

\subsubsection*{Feller's Proof of the `if' part of Classical
M\"{u}ntz Theorem.} Taking into account the decomposition
properties of signed measures, we see that in order to prove
$(b)\Rightarrow (a)$ in Theorem \ref{muntzmeasures} (which
corresponds to the `if' part of the classical M\"{u}ntz Theorem)
it suffices to prove the assertion for nonnegative measures $\mu$.
Now, we note that functions $f:(0,\infty)\to\mathbb{R}$ admitting
an expression of the form $(\ref{fmu})$ for a certain nonnegative
measure $\mu$ are completely monotone on $(0,\infty)$. This means
that they satisfy the inequalities
\[
(-1)^nf^{(n)}(t)\geq 0 \text{ for all } t>0 \text{ and all }
n=0,1,2,\ldots.
\]
(Indeed, it is a well known result by S.N. Bernstein
\cite{bernsteinmonotone} that $f$ being completely monotone on
$(0,\infty)$ and of the form $(\ref{fmu})$ for a certain
nonnegative measure $\mu$ are equivalent claims). Let us now
assume that (a) is not true, and suppose first that
$(\lambda_k)\uparrow \infty$. Under these conditions we can use
the following theorem:

\begin{theorem}[Feller, 1968]\label{teoremafeller} Let us assume that
$0<\lambda_0<\lambda_1<\lambda_2<\cdots $ with
$\lambda_n\to\infty$, and $\sum_{n=0}^{\infty}1/\lambda_n=
\infty$, and let $f:(0,\infty)\to\mathbb{R}$ be a completely
monotone function. Then
\begin{equation}\label{serienewtonmonotona}
f(t)=\sum_{n=0}^{\infty}f[\lambda_0,\lambda_1,\ldots,\lambda_n](t-\lambda_0)(t-\lambda_1)\cdots(t-\lambda_{n-1}),
\end{equation}
where the series is absolutely convergent for all $t>0$.
\end{theorem}

This proves that if $f(\lambda_k)=0$ for all $k$, then $f(k)=0$
for all $k$, and this means that the integral of any polynomial
against $\mu$ is zero, hence $\mu$ is zero, so (b) is also false.

Finally, we can use Morera's theorem to prove that $f(z)=
\int_{0}^{1}t^{z}\dd \mu (t)$ is holomorphic on the half plane
$\{z:{\rm Re} z>0\}$, so that the well known principle of identity
shows that $f$ is completely determined by its values on any
increasing bounded sequence $(\lambda_k)_{k=0}^{\infty}$. This
ends the proof of $(b)\Rightarrow (a)$ in Theorem
\ref{muntzmeasures}.

\Proof of Theorem \ref{teoremafeller}. It follows from the fact
that $f$ is completely monotone and
(\ref{diferenciasdivididasderivadas})  that
\[
(-1)^nf[\lambda_0,\lambda_1,\ldots,\lambda_n]\geq 0 \text{ for all
} n\geq 0.
\]
Let us now assume that $t\in [0,\lambda_0)$. Then all terms of the
series
\begin{equation}\label{seriemonotona}
\sum_{n=0}^{\infty}f[\lambda_0,\lambda_1,\ldots,\lambda_n](t-\lambda_0)(t-\lambda_1)\cdots(t-\lambda_{n-1})
\end{equation}
are positive. Setting
\[P_n(t):=\sum_{k=0}^{n}f[\lambda_0,\lambda_1,\ldots,\lambda_k](t-\lambda_0)(t-\lambda_1)\cdots(t-\lambda_{k-1})\]
and $R_n(t):=f(t)-P_n(t)$, we obtain
\[
R_n(t)=f[t,\lambda_0,\lambda_1,\ldots,\lambda_n](t-\lambda_0)(t-\lambda_1)\cdots(t-\lambda_n),
\]
so that
\begin{equation}\label{pcrece}
P_0(t)\leq P_1(t)\leq \cdots \leq  f(t)
\end{equation}
and there exists a function $\alpha(t)$ such that
\[
R_n(t)\downarrow \alpha(t).
\]
We want to show that $\alpha(t)=0$. Now, for $0<s<\lambda_0$, we
have that
\[
0\leq R_n(s)\leq
f[s,\lambda_0,\lambda_1,\ldots,\lambda_n](-1)^{n+1}\lambda_0\lambda_1\cdots\lambda_{n},
\]
so that
\begin{equation}\label{pruebamonotona}
\frac{s}{\lambda_n}\alpha(s)\leq \frac{s}{\lambda_n}R_n(s)\leq
f[s,\lambda_0,\lambda_1,\ldots,\lambda_n](-1)^{n+1}s\lambda_0\lambda_1\cdots\lambda_{n-1}.
\end{equation}
Now, the right side of (\ref{pruebamonotona}) is the $n$th term of
the series (\ref{serienewtonmonotona}) evaluated at $t=0$ when the
point $s$ is added to the sequence
$(\lambda_k)_{k=0}^{\infty}$. It follows that the series
$\sum_{k=0}^{\infty}\alpha(s)/\lambda_k$ is convergent, which is
consistent with our hypotheses on the sequence
$(\lambda_k)_{k=0}^{\infty}$ only if $\alpha(s)=0$. This proves
(\ref{serienewtonmonotona}) for all $t\in (0,\lambda_0)$. The same
argument works for $t\in (\lambda_{2k-1},\lambda_{2k})$ except
that the inequalities (\ref{pcrece}) can be asserted only for
$n\geq 2k$. On the intervals $(\lambda_{2k},\lambda_{2k+1})$ the
inequalities are reversed. This ends the proof, since
$\lim_{k\to\infty}\lambda_k= \infty$ implies that all points $t>0$
have been already considered.\endproof

\subsubsection*{Hirschman-Widder's and Gelfond's proof of the
`if' part of the M\"{u}ntz Theorem} The most famous proof of the
Weierstrass Approximation Theorem is based on the use of the
Bernstein polynomials:
\[
B_nf(x):=\sum_{k=0}^nf\!\left(\frac{k}{n}\right)\binom{n}{k}x^k(1-x)^{n-k}.
\]
Thus, it was an interesting (and difficult!) problem to find out
whether a suitable generalization of the Bernstein polynomials
would give a new proof of the M\"{u}ntz Theorem. This question was
solved in the positive by Hirschman and Widder \cite{widder} in
1949. Moreover, their proof was modified and extended by
A.~O.~Gelfond \cite{gelfond} in 1958 and included by G.~G.~Lorentz
in his book on Bernstein polynomials \cite[pp.
46--47]{lorentzbernstein}. In this subsubsection, we follow the
discussion in Lorentz's monograph.

The polynomials that will play the role of Bernstein polynomials
are defined in terms of the sequence of exponents
$(\lambda_k)_{k=0}^{\infty}$ as follows: Given $n,k\in\mathbb{N}$
such that $k\leq n$, we set
\begin{equation}\label{gelfondbasis}
g_{n,k}(x):=(-1)^{n-k}\lambda_{k+1}\cdots\lambda_n\sum_{i=k}^n\frac{x^{\lambda_i}}{(\lambda_i-\lambda_k)\cdots
(\lambda_i-\lambda_{i-1})(\lambda_i-\lambda_{i+1})\cdots
(\lambda_i-\lambda_{n})}
\end{equation}
and, given $f\in C[0,1]$, we set
\[
\eta_{n,k}:=\left[(1-\frac{\lambda_1}{\lambda_{k+1}})\cdots
(1-\frac{\lambda_1}{\lambda_{n}})\right]^{\frac{1}{\lambda_1}}
\quad \text{ for } \ 0\leq k<n,\  \text{ and } \ \eta_{n,n}:=1,
\]
and
\begin{equation}\label{gelfondpolynomials}
B_n^{\Lambda}(f)(x):=\sum_{k=0}^nf(\eta_{n,k})g_{n,k}(x).
\end{equation}
We can now state and prove the main result:
\begin{theorem}[Hirschman-Widder \cite{widder}, and Gelfond \cite{gelfond}]\label{teogelfond} Let
$f\in C[0,1]$ and assume that
\[0<\lambda_1<\lambda_2<\cdots,\quad
\lim_{k\to\infty}\lambda_k=\infty\quad \mbox{and}\quad
\sum_{k=1}^{\infty}\frac{1}{\lambda_k}=\infty.\]
 Then
\[
\lim_{n\to\infty}\|B_n^{\Lambda}(f)-f\|_{[0,1]}=0.
\]
\end{theorem}
\Proof. Let us consider, for  the function $f(z)=x^z$, its divided
differences with respect to the nodes
$(\lambda_k)_{k=1}^{\infty}$. It is clear that
\begin{equation}\label{gelfonddivididas}
g_{n,k}(x)=(-1)^{n-k}\lambda_{k+1}\cdots\lambda_nf[\lambda_k,\lambda_{k+1},\ldots,\lambda_n].
\end{equation}
In particular, this implies that
\[
g_{n,k}(x)=(-1)^{n-k}\lambda_{k+1}\cdots\lambda_n \frac{1}{2\pi
\ii}\int_{C}\frac{x^z\dd z}{(z-\lambda_k)\cdots(z-\lambda_n)},
\]
where $C$ is any simple closed curve that contains the nodes
$(\lambda_i)_{i=k}^n$ in its interior Int$(C)$, and such that
$f(z)=x^z$ is holomorphic in a neighborhood of ${\rm Int} (C)\cup
C$. Now we prove a few technical results:

\begin{lemma} The polynomials $\{g_{n,k}\}_{k=0}^n$ form a
partition of unity on $[0,1]$.
\end{lemma}
\Proof. Taking into account (\ref{gelfonddivididas}) and
(\ref{diferenciasdivididasderivadas}), we get
\[
g_{n,k}(x)=\frac{\lambda_{k+1}\cdots
\lambda_n}{(n-k)!}x^{\tau}(-\log x)^{n-k}\geq 0.
\]
Moreover, taking into account the identity (easily checked by
induction on $n$)
\[
\frac{1}{z}=\frac{1}{z-\lambda_n}-\frac{\lambda_n}{
(z-\lambda_{n-1})(z-\lambda_n)}+\cdots+(-1)^{n}\frac{\lambda_1
\cdots \lambda_n}{z(z-\lambda_1)\cdots (z-\lambda_n)}
\]
and multiplying by ${x^z}/(2\pi \ii)$ and integrating along $C$,
we get
\[
1=\frac{1}{2\pi \ii}\int_{C}\frac{x^z}{z}\dd
z=\sum_{k=0}^ng_{n,k}(x),
\]
which is what we wanted to prove.\endproof

\begin{lemma}The following identities hold:
\begin{equation}\label{xlambda1}
x^{\lambda_1}=\sum_{k=0}^n\eta_{n,k}^{\lambda_1}g_{n,k}(x)
\end{equation}
and
\begin{equation}\label{xdoslambda1}
x^{2\lambda_1}=\sum_{k=1}^n\eta_{n,k}^*g_{n,k}+\frac{1}{2}\eta_{n,0}^*g_n^*(x),
\end{equation}
where
\[
\eta_{n,k}^*:=(1-\frac{2\lambda_1}{\lambda_{k+1}})\cdots(1-\frac{2\lambda_1}{\lambda_n})
\]
and $g_n^*$ is the polynomial $g_{n+1,1}$ associated with the
nodes
$\lambda_0=0,\lambda_1,2\lambda_1,\lambda_2,\ldots,\lambda_n$
$($taken in increasing order$)$.
\end{lemma}
\Proof. The idea is analogous to that in the previous
lemma. We write $1/(z-\lambda_1)$ in a different way (this is
again easy to check by induction on $n$):
\[
\frac{1}{z-\lambda_1}=\frac{1}{z-\lambda_n}-\frac{\lambda_n-\lambda_1}{
(z-\lambda_{n-1})(z-\lambda_n)}+\cdots+(-1)^{n-1}\frac{(\lambda_2-\lambda_1)
\cdots (\lambda_n-\lambda_1)}{(z-\lambda_1)\cdots (z-\lambda_n)}
\]
and multiplying by $x^z/(2\pi\ii)$ and integrating over $C$, we
get (\ref{xlambda1}). To prove (\ref{xdoslambda1}), we use the
same arguments but based on the formula
\[
\frac{1}{z-2\lambda_1}=\frac{1}{z-\lambda_n}-\frac{\lambda_n-2\lambda_1}{
(z-\lambda_{n-1})(z-\lambda_n)}+\cdots+(-1)^{n}\frac{\lambda_1(\lambda_2-2\lambda_1)
\cdots
(\lambda_n-2\lambda_1)}{(z-\lambda_1)(z-2\lambda_1)(z-\lambda_2)\cdots
(z-\lambda_n)}.
\]
\endproof

\noindent Let us continue with the proof of Theorem
\ref{teogelfond}. Consider the functions $T_n$ given by
\[
T_n(x):=
\sum_{k=0}^n(x^{\lambda_1}-\eta_{n,k}^{\lambda_1})^2g_{n,k}(x).
\]
Then
\begin{eqnarray*}
T_n(x) &=&
x^{2\lambda_1}-2x^{\lambda_1}\sum_{k=0}^n\eta_{n,k}^{\lambda_1}g_{n,k}(x)+
\sum_{k=0}^n\eta_{n,k}^{2\lambda_1}g_{n,k}\\
&=&
\sum_{k=1}^n(\eta_{n,k}^{2\lambda_1}-\eta_{n,k}^*)g_{n,k}(x)-\frac{1}{2}\eta_{n,0}^*g_n^*(x).
\end{eqnarray*}
Now, the sequence $(\eta_{n,0}^*)$ converges to zero since the
product $\prod_{n=1}^{\infty}(1-2\lambda_1/\lambda_n)$ diverges to
zero because of our hypothesis on the sequence
$(\lambda_k)_{k=1}^{\infty}$. This means that
$-\frac{1}{2}\eta_{n,0}^*g_n^*(x)$ converges uniformly to zero,
since we know that $0\leq g_n^*(x)\leq 1$ for all $x\in [0,1]$.

Let us now show that
\begin{equation}\label{cosa}
\eta_{n,k}^{2\lambda_1}-\eta_{n,k}^*\to 0
\end{equation}
uniformly in $k\geq 1$. To prove this, let $\varepsilon>0$ be
arbitrary and fix $n_0$ such that, for all $k\geq n_0$,
$2\lambda_1<\lambda_k$ and
\[
(1+\varepsilon)\log(1-\frac{\lambda_1}{\lambda_k})^2\leq
\log(1-\frac{2\lambda_1}{\lambda_k})\leq
\log(1-\frac{\lambda_1}{\lambda_k})^2.
\]
(This is possible since $\lambda_k\to\infty$, $\log$ is an
increasing function and $\log t <0$ for $t\in (0,1)$.) Then
\begin{equation}\label{cosa2}
(1-\frac{\lambda_1}{\lambda_k})^{2(1+\varepsilon)}\leq
(1-\frac{2\lambda_1}{\lambda_k})\leq
(1-\frac{\lambda_1}{\lambda_k})^2.
\end{equation}
Since the products
$\prod_{k=1}^\infty(1-\frac{\lambda_1}{\lambda_k})$ and
$\prod_{k=1}^\infty(1-\frac{2\lambda_1}{\lambda_k})$ are both
divergent to zero, we can consider only the values $k\geq n_0$ in
$(\ref{cosa})$, so that we can use $(\ref{cosa2})$ for the factors
representing $\eta_{n,k}^{2\lambda_1}$ and $\eta_{n,k}^{*}$. If
$\eta_{n,k}^{2\lambda_1}\leq\varepsilon $ then also
$\eta_{n,k}^{*}\leq\varepsilon $ and $|\eta_{n,k}^{2\lambda_1}-
\eta_{n,k}^{*}|\leq 2\varepsilon$. On the other hand, if
$\eta_{n,k}^{2\lambda_1}\geq\varepsilon $ then
\[
0\leq \eta_{n,k}^{2\lambda_1}- \eta_{n,k}^{*}\leq
\eta_{n,k}^{2\lambda_1}-
(\eta_{n,k}^{2\lambda_1})^{1+\varepsilon}\leq
1-(\eta_{n,k}^{2\lambda_1})^{\varepsilon}\leq
1-\varepsilon^{\varepsilon},
\]
which is arbitrarily small for $\varepsilon\to 0$. This proves
$(\ref{cosa})$ and, consequently, that $T_n(x)\to 0$ uniformly in
$x\in [0,1]$.

Let us take $f\in C[0,1]$ and let $M\geq \|f\|_{[0,1]}$.
Obviously, $f$ is uniformly continuous so that we may assume that
for a given $\varepsilon>0$ we have chosen $\delta>0$ such that
$|x-y|<\delta$ implies
$|f(x^{1/\lambda_1})-f(y^{1/\lambda_1})|<\varepsilon$. Then
\[
B_n^{\Lambda}(f)(x)-f(x)=\sum_{k=0}^n(f(\eta_{n,k})-f(x))g_{n,k}(x)
\]
(recall that the polynomials $\{g_{n,k}\}_{k=0}^n$ are a partition
of unity in $[0,1]$). Hence we can decompose the summation formula
into two parts: the first one containing those indices $k$ such
that $|\eta_{n,k}^{\lambda_1}-x^{\lambda_1}|\leq \delta$ and the
second one, where this inequality does not hold. The first part of
the summation formula is $\leq \varepsilon$ and for the second
part we use that $|f(\eta_{n,k})-f(x)|\leq 2M$ and
$(\eta_{n,k}^{\lambda_1}-x^{\lambda_1})^2/\delta^2\geq 1$ to
conclude that $2MT_n(x)/\delta^2$ is an upper bound of this part.
This obviously implies that $B_n^{\Lambda}(f)(x)$ converges to
$f(x)$ uniformly on $[0,1]$.\endproof

\subsection{Proof of M\"{u}ntz theorem via complex analysis}
In this subsection, we will use some basic results from complex
analysis to give another proof of the Müntz Theorem. It is because
of our use of a complex variable that we introduce a minor
modification in the space of functions we want to approximate.
Concretely, we assume that our space of functions is
$C([0,1],\mathbb{C})$ and M\"{u}ntz polynomials have complex
coefficients. Clearly, the Müntz Theorem corresponding to this
context is the following one (which is equivalent to the classical
Müntz Theorem since the variable $z$ runs on $[0,1]$ which is a
subset of $\mathbb{R}$, so that in order to approximate a
continuous function $f(z)=u(z)+{\ii}v(z)$ with a complex
polynomial $p(z)=\sum_{i=0}^n\alpha_iz^{\lambda_i}$ we only need
to choose the coefficients $\alpha_i=a_i+{\ii}b_i$ in such a form
that $\sum_{i=0}^na_iz^{\lambda_i}$ approximates $u(z)$ and
$\sum_{i=0}^nb_iz^{\lambda_i} $ approximates $v(z)$):

\begin{theorem}[M\"{u}ntz Theorem for Complex-Valued Functions]\label{muntzcomplex}
Let\textbf{\ }$\Lambda =(\lambda _{k})_{k=0}^{\infty }$,\ \
$0=\lambda _{0}<\lambda _{1}<\cdots $ be an increasing sequence of
non-negative real numbers. Then $\Pi_{\mathbb{C}}(\Lambda)$ is a
dense subset of $C([0,1],\mathbb{C})$ if and only if
$\sum_{k=1}^{\infty}1/\lambda _{k}=\infty$.
\end{theorem}

Let $\mathbb{D}:=\{z\in \Bbb{C}:|z|<1\}$ be the unit disc, and
\[
H^{\infty }(\mathbb{D}):=\{f:f\text{ is holomorphic on
}\mathbb{D}\text{ and }\|f\|_{H^{\infty }(\mathbb{D})}=\sup_{z\in
\mathbb{D}}|f(z)|<\infty \}
\]
be the algebra of bounded analytic functions defined on
$\mathbb{D}$. The proof we present of Theorem \ref{muntzcomplex}
(which is due to Feinerman and Newman \cite{feinerman}) is based
on the following lemmas.

\begin{lemma}[Blaschke Products]\label{blaschkeproduct} The function $f:\Bbb{D}\rightarrow \Bbb{C}$ \ belongs
to $H^{\infty }(\mathbb{D})$ if and only if it can be decomposed
as
\[
f(z)=z^{p}\prod_{k=0}^{\infty }(z-\lambda _{k})h(z)
\]
for a certain choice of a natural number $p\geq 0$, a sequence
of points $(\lambda _{k})\subset \mathbb{D}$ such that
$\sum_{k=0}^{\infty }(1-|\lambda _{k}|)<\infty $, and a function
$h\in H^{\infty }(\mathbb{D})$ without zeros on $\mathbb{D}$.
\end{lemma}

\Proof. See \cite[pages 63--67]{hoffman}.
\endproof

\begin{lemma}\label{lemacomplejo}
If\/ $\sum_{k=0}^\infty  1/\lambda_k=\infty $ and $\eta $ is a
complex Borel measure on $[0,1]$ such that
\[
\int_{0}^{1}t^{\lambda _{k}}\dd \eta (t)=0, \qquad k=0,1,\ldots,
\]
then
\[
\int_{0}^{1}t^{k}\dd \eta (t)=0, \qquad k=0,1,\ldots\,.
\]
\end{lemma}

\Proof. We may assume, without loss of generality, that our
measure is concentrated on $(0,1]$. Then we use Morera's theorem
to prove that $h(z):=\int_{0}^{1}t^{z}\dd\eta (t)$
is holomorphic on the half plane $\{z:\text{Re} z>0\}$ and $%
h(\lambda _{k})=0$, $k=0,1,\ldots$. Moreover, $h$ is bounded on
the half plane,
since if we decompose $z=x+\ii y$, then $|t^{z}|=t^{x}\leq 1$ for all $%
t\in \lbrack 0,1]$. It follows that $g(z):=h((1+z)/(1-z))\in H%
^{\infty }(\mathbb{D})$ and $g(\alpha _{k})=0$, $k=0,1,\ldots$,
where $\alpha _{k}:=(\lambda _{k}-1)/(\lambda _{k}+1)$, for all
$k$.

Now, it is clear that $\sum 1/\lambda _{k}=\infty $ (which is our
hypothesis), implies that $\sum (1-|\alpha _{k}|)=\infty $ , so
that $g=0$. This of course implies that $h(k)=0$ for all $k\in
\mathbb{N}$.
\endproof

\Proof{Theorem \ref{muntzcomplex}}. Let $Y:=\overline{\Pi (\Lambda
)}$ be the closure of $\Pi (\Lambda )$
in $X = C[0,1].$ It follows from the Hahn-Banach Theorem that $%
f\in X\setminus Y$ if and only if there exists a bounded linear
functional $L\in X^*$ such that
$L|_{Y}=0$ but $L f\neq 0$. Now, if $L\in C^*[0,1]$ then $%
Lf=\int_{0}^{1}f(t)\dd \eta $ for all $f$ and a certain finite
complex Borel measure $\eta $ defined on $[0,1]$. It follows from
Lemma \ref{lemacomplejo} (and the Weierstrass Approximation
Theorem) that $Lf=0$ for all $f$, which is in contradiction to our
hypotheses. This proves that $\sum 1/\lambda _{k}=\infty $ is a
sufficient condition for the density of $\Pi (\Lambda )$ in
$C[0,1]$.

We now prove that the condition is also necessary. Let us assume
that $\sum  1/\lambda_k <\infty $, and let us define the function:
\[
f(z):=\frac{z}{(2+z)^{3}}\prod_{k=1}^{\infty }\frac{\lambda
_{k}-z}{2+\lambda _{k}+z}\,.
\]
Now,
\[
1-\frac{\lambda _{k}-z}{2+\lambda _{k}+z}=\frac{2+2z}{2+\lambda
_{k}+z},
\]
so that the infinite product that appears in the definition of $f$
converges uniformly on compact subsets of $\Bbb{C}\setminus
(\{-2\}\cup \{-2-\lambda _{k}\}_{k=0}^{\infty })$. Hence $f$ is a
meromorphic function on $\Bbb{C}$ with poles $\{-2\}\cup
\{-2-\lambda _{k}\}_{k=0}^{\infty }$ and zeros $\{0\}\cup
\{\lambda _{k}\}_{k=0}^{\infty }$. Furthermore, each
factor of our infinite product is (in absolute value) less than $1$ for all $%
z$ such that $\text{Re}\, z>-1$. On the other hand, the
restriction of $f $ to the line $\text{Re}\, z=-1$ is an
absolutely integrable function (this follows
from the fact that we have divided by $(2+z)^{3}$). Let us fix $z$ with $\text{Re}%
\, z>-1$ and consider the Cauchy formula for $f$ taking as path of
integration the circle centered at $-1$ of radius $R>1+|z|$, from $-1-%
\ii R$ to $-1+\ii R$, plus the interval $[-1+\ii R,-1-%
\ii R]$. If we let $R\rightarrow \infty $ \ then we can eliminate
the part of the formula associated with the semicircle (note that
$|2+z|^3>R^3$ there), and we obtain
\begin{eqnarray}
f(z) &=&-\frac{1}{2\pi }\int_{-\infty }^{+\infty }\frac{f(-1+ \ii s)\dd s}{%
-1+\ii s-z}  \label{++} \\[10pt]
&=&\int_{0}^{1}t^{z}\left\{ \frac{1}{2\pi }\int_{-\infty }^{+\infty }f(-1+%
\ii s)\exp (-\ii s\log t)\dd s\right\} \dd t,\text{ } \nonumber
\end{eqnarray}
since
\[
\frac{1}{-1+ \ii s-z}=\int_{0}^{1}t^{z- \ii s}\dd t=\int_{0}^{1}t^{z}%
\exp (-\ii s\log t)\dd t.
\]
If we define $g(s):=$ $f(-1+\ii s)$ then the inner integral which
appears in formula (\ref{++}) is $\widehat{g}(\log t)$, where
$\widehat{g}$ denotes the Fourier transform of $g$ which is
clearly a continuous bounded function defined on $(0,1]$, so that
if $\dd \eta (t)=\widehat{g}(\log t)\dd t$ then $\eta$ is a
 complex Borel measure. Therefore, the
bounded linear functional $h\mapsto\int h\dd \eta$ annihilates $Y$
but it is not identically zero, hence $Y\neq X$ whenever
$\sum_{k=0}^\infty 1/\lambda_k <\infty $, which is what we wanted
to prove.\endproof

\section{The Full M\"{u}ntz theorem and polynomial inequalities}

The Classical M\"{u}ntz Theorem was only stated for increasing
sequences of nonnegative real numbers $0=\lambda _{0}<\lambda_{1}<
\cdots$.  It would be interesting to know if a general result,
dealing with arbitrary sequences of exponents, is possible. We
call such a result a \dword{Full M\"{u}ntz Theorem}. Moreover, it
would be interesting to know such a result not only for the space
$C[0,1]$ but also for $C(K)$ with $K$ a compact subset of
$\mathbb{R}$ or $\mathbb{C}$ and for other function spaces such as
$L^p[a,b]$, etc. As we have already noted, Sz\'{a}sz' proof of the
M\"{u}ntz Theorem in 1916 included sequences of exponents more
general than those treated by M\"{u}ntz in his paper of 1914. In
particular, for real exponents, he proved a result that works
whenever $\lim\inf_{k\to\infty}\lambda_k>0$, so that Sz\'{a}sz'
theorem was not properly a Full M\"{u}ntz theorem although it was
near to it. It is remarkable that, after Sz\'{a}sz' work, the
search for a Full M\"{u}ntz theorem was not immediate. Moreover,
since the appearance of the M\"{u}ntz and Sz\'{a}sz theorems it
was clear that the role of the origin in these results was very
important and to extend the results to spaces of functions defined
away from the origin was a nontrivial task. Moreover, nobody had
answered the main question of characterizing the elements in the
closure of a M\"{u}ntz space $\Pi(\Lambda)$ when this space is not
a dense subset of $C[a,b]$.

It was Fields medalist Laurent Schwartz \cite{schwartz} who proved
a Full M\"{u}ntz theorem for the space of square-integrable
functions $L^2[a,b]$ on general intervals $[a,b]$. Schwartz also
characterized the density of M\"{u}ntz spaces in $C[a,b]$ for
$0\not\in [a,b]$ and conjectured a necessary and sufficient
condition for the density in $C[0,1]$ of $\Pi(\Lambda)$ with
general sequences, but he did not prove the result.

It was only a few years later when Siegel \cite{siegel}, in a
beautiful paper where he included a difficult generalization of
Sz\'{a}sz' theorem, proved Schwartz' conjecture for the first time
using complex variable techniques. The deepest work related to the
Full M\"{u}ntz Theorem has only recently been done by P.~Borwein
and T.~Erd\'elyi, sometimes in collaboration with several other
authors. They proved that there is a strong connection between
density results for M\"{u}ntz spaces and the study of some special
inequalities for these polynomials. In fact, their book
``Polynomials and Polynomial Inequalities'' \cite{borweinerdelyi}
constitutes a deep contribution to this subject and contains a
guided investigation of the Full M\"{u}ntz Theorem. The reader
should attempt the solution of the exercises in Chapter 4 of that
book or, avoiding much effort, read this section where we will
concentrate almost all our attention on the study of the Full
M\"{u}ntz Theorem for the spaces $C[0,1]$, $L^p[0,1]$ and $C[a,b]$
(with $0<a<b$). We also include some ideas related to the proof by
Borwein and Erd\'elyi of a Full M\"{u}ntz Theorem for the space
$C(K)$, where $K$ is a compact set with positive Lebesgue measure,
and a recent result by the author where a Full M\"{u}ntz Theorem
is proved for the space of continuous functions on quite general
countable compact sets.

\subsection{Full M\"{u}ntz theorem on $[0,1]$}
We start with the precise statement of the main results of this
section.

\begin{theorem} {\bf (Full M\"{u}ntz Theorem for $C[0,1]$)} \label{muntzc01} Let us
assume that $\Lambda =(\lambda _{k})_{k=1}^{\infty }$ is a
sequence of distinct real positive numbers. Then $\Pi (\Lambda
\cup \{0\})$ is dense in $C[0,1]$ if and only if
\begin{equation}\label{condmuntzc01}
\sum_{k=1}^{\infty }\frac{%
\lambda _{k}}{\lambda _{k}^{2}+1}=\infty.
\end{equation}
\end{theorem}

\begin{theorem} \label{muntzl201}{\bf (Full M\"{u}ntz Theorem for $L_2[0,1]$)} Let us assume that $\Lambda
=(\lambda _{k})_{k=1}^{\infty }$ is a sequence of distinct real
numbers greater than $-1/2$. Then $\Pi (\Lambda )$ is dense in
$L_{2}[0,1]$ if and only if
\begin{equation}\label{condmuntzl201}
\sum_{k=1}^{\infty }%
\frac{2\lambda_{k}+1}{(2\lambda_{k}+1)^{2}+1}=\infty.
\end{equation}
\end{theorem}

Although we will not prove it in this paper, we would like to
include here the corresponding result for $L^p[0,1]$:

\begin{theorem} \label{muntzlp01} {\bf (Full M\"{u}ntz Theorem for $L_p[0,1]$)} Let
$p\in (0,\infty )$ and let us assume that $\Lambda
=(\lambda_{k})_{k=1}^{\infty }$ is a sequence of distinct real
numbers greater than $-1/p$. Then $\Pi (\Lambda )$ is dense in
$L_{p}[0,1]$ if and only if
\begin{equation}\label{condmuntzlp01}
\sum_{k=1}^{\infty }\frac{\lambda_{k}+1/p}{(\lambda _{k}+1/p)^{2}+1%
}=\infty.
\end{equation}
\end{theorem}

Theorem \ref{muntzc01} was conjectured by Schwartz, proved by
A.~R.~Siegel \cite{siegel} for the first time and, following the
ideas introduced by Sz\'asz in his famous 1916 paper, reproved by
Borwein and Erd\'elyi. Theorem \ref{muntzl201} was proved by
Sz\'asz \cite{szasz}. Theorem \ref{muntzlp01} was proved for $p=1$
and conjectured for $p>0$ by Borwein and Erd\'elyi
\cite{borweinerdelyi}. Moreover, it was proved by Operstein
\cite{operstein} for the case $1<p<\infty$, by Erd\'elyi and
Johnson \cite{erdejohn} (using quasi-Banach space theory) for
$0<p<\infty$ and, quite recently, by Erd\'elyi
\cite{erdelyiconstructiv05} for $0<p<\infty$ with an
``elementary'' proof.

In this section, we will prove Theorem \ref{muntzc01}. As a first
step, we prove the easier Theorem \ref{muntzl201} and we use it to
prove some particular cases of Theorem \ref{muntzc01}. Then we
introduce several polynomial inequalities that will be needed to
complete the proof of this theorem. These inequalities are proved
in the second part of this section. The third and fourth part are
devoted to a characterization of the closure of nondense M\"{u}ntz
subspaces of $C[0,1]$ and the statement and proof of the Full
M\"{u}ntz Theorem for intervals $[a,b]$ away from the origin.

Finally, the fifth and sixth part are devoted to the Full
M\"{u}ntz theorem for $C(K)$ for  compact sets $K\subset
[0,\infty)$ more general than intervals.


\Proof {Theorem \ref{muntzl201}}. We have already proved the
formula:
\[
E(x^{q},\Pi (\Lambda _{n}))_{L^{2}(0,1)}=\frac{1}{\sqrt{2q+1}}\prod_{k=1}^{n}%
\frac{\left| q-\lambda _{k}\right| }{|q+\lambda _{k}+1|}.
\]
Hence $x^{q}\in \overline{\Pi (\Lambda )}^{L^{2}(0,1)}$ if and
only if
\[
\lim_{n\rightarrow \infty }\prod_{k=1}^{n}\left| \frac{q-\lambda _{k}}{%
q+\lambda _{k}+1}\right| \text{ }\left( =\lim_{n\rightarrow \infty
}\prod_{k=1}^{n}\left| 1-\frac{2q+1}{q+\lambda _{k}+1}\right|
\right) =0.
\]
We decompose the above product making a distinction between the cases $%
\lambda _{k}\in (-1/2,q]$ and $\lambda _{k}\in (q,\infty )$, which
leads us to the following reformulation of the above condition:
\[
\lim_{n\rightarrow \infty }\prod_{k\leq n;\lambda _{k}\in (-1/2,q]\text{ }%
}\left| 1-\frac{2q+1}{q+\lambda _{k}+1}\right|\ \  \prod_{k\leq n;\text{ }%
\lambda _{k}\in (q,\infty )}\left| 1-\frac{2q+1}{q+\lambda
_{k}+1}\right| =0,
\]
which is clearly equivalent to stating that
\[
\sum_{k\geq 1\text{; }\lambda_{k}\in (q,\infty )}\frac{1}{2\lambda_{k}+1}%
=\infty \quad\text{ \ or } \quad\sum_{k\geq 1\text{; }\lambda
_{k}\in (-1/2,q]}(2\lambda_{k}+1)=\infty,
\]
and this can be rewritten as
\[
\sum_{k=1}^{\infty }\frac{2\lambda _{k}+1}{(2\lambda
_{k}+1)^{2}+1}=\infty,
\]
which is what we wanted to prove.\endproof

We subdivide the proof of Theorem \ref{muntzc01} into several
cases, depending on which of the following three conditions are
satisfied by the sequence $(\lambda_{k})$:

\begin{itemize}
\item[\textbf{H1}]  \textbf{\ }$\inf_{k\in
\mathbb{N}}\lambda_{k}>0$.

\item[\textbf{H2}]  $\lim_{k\rightarrow \infty }\lambda_{k}=0.$
(In this case, the identities $\sum \frac{\lambda_{k}}{\lambda
_{k}^{2}+1} =\infty $ and $\sum \lambda_{k}=\infty $ are
equivalent.)

\item[\textbf{H3}]  $(\lambda_{k})=(\alpha_{k})\cup (\beta _{k})$,
with $\alpha _{k}\rightarrow 0$ and $\beta _{k}\rightarrow \infty
$.
(In this case, the identities $\sum \frac{\lambda_{k}}{%
\lambda_{k}^{2}+1}=\infty $ and $\sum \alpha_{k}+\sum
\frac{1}{\beta_{k}}=\infty $ are equivalent.)
\end{itemize}

\noindent%
%

At first glance it seems that condition H3 is not significant
since, taking subsequences if necessary, H1 and H2 produce all
possibilities. Theorem \ref{muntzc01} has been stated in a
unified form that depends on the convergence character of the series $%
\sum \frac{\lambda_{k}}{\lambda_{k}^{2}+1}$, and
this depends on the knowledge of the boundedness character of the set $%
\{\lambda_{k}\}_{k=0}^{\infty }$. In particular, for the case
$\inf_{k\geq 0}\lambda_k=0$ and in order to characterize for which
sequences we have $\sum
\frac{\lambda_{k}}{\lambda_{k}^{2}+1}=\infty$, we must take into
account both possibilities H2 and H3. We will see that the study
of the case described by H3 is precisely the most difficult to
handle.

In fact, H1 -- H3 do not cover all cases (e.g. if there is a
subsequence converging to $0$ and another one converging to $1$)
but in the missing cases there is a subsequence that is bounded
away from $0$ and infinity, and in this case (\ref{condmuntzc01})
is automatically true, and so is the denseness by H1 applied to
this subsequence.

\Proof {Theorem \ref{muntzc01}, using Theorem \ref{muntzl201}, and
assuming that $\inf_{k\in \mathbb{N}}\protect\lambda _{k}>0$}. Let
$0<\delta \leq $ $\inf_{k\in \mathbb{N}}\lambda_{k}$. We make the
change of variable $x\rightarrow x^{\frac{1}{\delta }}$ and solve
the problem for exponents $\lambda_{k}^{\ast }=\lambda_{k}/\delta
$ that satisfy $\inf_{k\in \mathbb{N}}\lambda_{k}^{\ast }\geq 1$.
This means that we may assume, without loss of generality, that
$\inf_{k\in \mathbb{N}}\lambda _{k}\geq 1$. Then
\[
\sum_{k=1}^{\infty }\frac{\lambda_{k}}{\lambda_{k}^{2}+1}=\infty
\qquad \text{ if and only if }\qquad \sum_{k=1}^{\infty
}\frac{(\lambda_{k}-1)}
{(\lambda_{k}-1)^{2}+1}=\text{ }\sum_{k=1}^{\infty }\frac{2(\lambda_{k}-1)+1}{%
(2(\lambda_{k}-1)+1)^{2}+1}=\infty.
\]
Let us first assume that
\[\sum \frac{2(\lambda_{k}-1)+1}{%
(2(\lambda_{k}-1)+1)^{2}+1}=\infty .\]
 It then follows from Theorem \ref{muntzl201} that $%
\Pi ((\lambda_{k}-1)_{k=1}^{\infty })$ is dense in $L^{2}[0,1]$. Let now $%
q\in \mathbb{N}$ be arbitrarily chosen. We can use the Sz\'asz
trick as described by the inequalities (\ref{trucoszasz}) (see 
Section 2 of this paper) to prove that
\[E(x^{q},\Pi
(\Lambda _{n}))_{C[0,1]}\leq qE(x^{q-1},\Pi (\Lambda _{n}^{\ast
}))_{L^{2}[0,1]},\]
 where $\Lambda _{n}=(\lambda
_{k})_{k=1}^{n}$ and $\Lambda _{n}^{\ast }=(\lambda
_{k}-1)_{k=1}^{n}$. This error goes to zero for all choices $q>0$,
which proves that
condition $\sum \lambda_{k}/(\lambda_{k}^{2}+1)%
=\infty $ is sufficient for the density of $\Pi (\Lambda )$ in $C%
[0,1]$.

On the other hand, if $\Pi ((\lambda_{k})_{k=1}^{\infty })$ is
dense in $C[0,1]$ then, taking into consideration that $C[0,1]$ is
dense in $L^{2}[0,1]$ and $\left\| \cdot \right\|
_{L^{2}[0,1]}\leq \left\| \cdot \right\|
_{C[0,1]}$, we have that $\Pi ((\lambda _{k})_{k=1}^{\infty })$ is
also dense in $L^{2}[0,1]$. Hence, using Theorem \ref{muntzl201},
$\sum \lambda _{k}/(\lambda_{k}^{2}+1)=\infty $.\endproof

\Proof {Theorem \ref{muntzc01} when $\lambda
_{k}\rightarrow 0$}.  We start by noting that in this case the identities $%
\sum \lambda_{k}/(\lambda_{k}^{2}+1)=\infty $ and $%
\sum \lambda_{k}=\infty $ are equivalent. It follows from the
Hahn-Banach and Riesz Representation Theorems that ${\rm span\,}\{1,x^{\lambda _{k}}%
\mathbf{\}}_{k=1}^{\infty }$ is dense in $C[0,1]$ if and only if $%
\{\lambda _{k}\}$ is not a subset of the set of zeros of any
nontrivial function of the form
\[f_{\mu }(z)=\int t^{z}\dd \mu (t)\]
 for some
finite Borel measure $\mu $. This condition is equivalent to
saying that $\{(\lambda_{k}-1)/(\lambda_{k}+1)\}$ is not the zero
set of any function of the form
\[g(z):=f_{\mu }((1+z)/(1-z))\in H^{\infty }(\mathbb{D}).\]

Now, if $\lambda_{k}\rightarrow 0$, then the equation $%
\sum \lambda_{k}=\infty $ implies that
\[
\sum_{k=1}^{\infty }\left( 1-\left| \frac{\lambda_{k}-1}{\lambda_{k}+1}%
\right| \right) =\infty \text{,}
\]
so that, using Lemma \ref{blaschkeproduct}, it is clear that
$\{(\lambda
_{k}-1)/(\lambda_{k}+1)\}$ is not the zero set of any $g(z)\in H%
^{\infty }(\mathbb{D})$. This means that $\sum \lambda_{k}=\infty$
is a sufficient condition (whenever $\lambda_{k}\rightarrow 0$)
for the density of ${\rm span\,}\{1,x^{\lambda _{1}},
x^{\lambda_{2}},\ldots \mathbf{\}}$ in $C[0,1]$.

In order to prove that condition (\ref{condmuntzc01}) is also
necessary, we need to introduce the following theorem:

\begin{theorem}[Newman's Inequality]\label{newmaninequality} Assume that $\Lambda =(\lambda
_{k})_{k=1}^{\infty }$ is a sequence of distinct positive real
numbers. Then the inequality
\[
\left\| xp^{\prime }(x)\right\| _{[0,1]}\leq 11\left(
\sum_{k=0}^{n}\lambda_{k}\right) \left\| p(x)\right\| _{[0,1]}
\]
holds for all $p\in \Pi (\Lambda _{n})$ and all $n\in \mathbb{N}$.
\end{theorem}

\noindent%
%
If $M:=\sum_{k=1}^{\infty }\lambda _{k}<\infty $,  then we have
that
\[
\left\| xp^{\prime }(x)\right\| _{[0,1]}\leq 11M\left\|
p(x)\right\| _{[0,1]}
\]
for all $p\in \Pi (\Lambda )$, which contradicts the density of
$\Pi (\Lambda )$ in $C[0,1]$. For let us assume
that $\Pi (\Lambda )$ is dense in $C[0,1]$. If, for example, we set $%
f(x)=(1-x)^{1/2}$ then for every natural number $m$ there exists a $%
p\in \Pi (\Lambda )$ such that $||p-f||\leq 1/m^{2}$. Hence
\begin{eqnarray*}
|p(1-1/m^{2})-p(1)| &\geq &|f(1-1/m^{2})-1/m^{2}-(f(1)+1/m^{2})| \\
&=&1/m-2/m^{2},
\end{eqnarray*}
and it then follows from the Mean Value Theorem that
\begin{eqnarray*}
|\xi p^{\prime }(\xi )| &=&\xi
\frac{|p(1-1/m^{2})-p(1)|}{1/m^{2}}\geq
(1-1/m^{2})\frac{1/m-2/m^{2}}{1/m^{2}} \\
&=&(1-1/m^{2})(m-2)\geq \frac{m-2}{2}
\end{eqnarray*}
for a certain $\xi \in (1-1/m^{2},1)$. This clearly is in
contradiction with
\[\left\| xp^{\prime }(x)\right\| _{[0,1]}\leq
11M\left\| p(x)\right\| _{[0,1]},\]
 since $m$ is arbitrary.\endproof

\Proof{Newman's Inequality}. We may assume, without loss of
generality, that $\sum_{k=0}^{n}\lambda_{k}=1$ since we may make
the change of variable
\[x\rightarrow
x^{1/( \sum_{k=0}^{n}\lambda_{k})}.\]
 Set $x=\ee^{-t}$. If $%
p(x)=\sum_{k=0}^{n}a_{k}x^{\lambda_{k}}$ and $q(t)=p(\ee^{-t})=%
\sum_{k=0}^{n}a_{k}\ee^{-\lambda_{k}t}$ then
\[
xp^{\prime }(x)=\sum_{k=0}^{n}\lambda_{k}a_{k}x^{\lambda
_{k}}=\sum_{k=0}^{n}\lambda_{k}a_{k}\ee^{-\lambda_{k}t}=q^{\prime }(t)\text{,%
}
\]
so that we have changed our problem to one of estimating the
uniform norm, on the interval $[0,\infty )$, of the derivatives of
functions of the form
\begin{equation}
\sum_{k=0}^{n}a_{k}\ee^{-\lambda_{k}t}  \label{expo}
\end{equation}
in terms of their uniform norms in the same interval.

Let
\[B(z):=\prod_{k=0}^{n}\frac{z-\lambda _{k}}{z+\lambda_{k}}\]
 and define
\[
T(t):=\frac{1}{2\pi \ii}\int_{\Gamma }\frac{\ee^{-zt}}{B(z)}\dd
z,\qquad \text{ where }\quad\Gamma :=\{z:|z-1|=1\}.
\]
It follows from the residue theorem that $T$ is of the form
(\ref{expo}). To prove Newman's inequality we first prove the
following estimate:
\begin{equation}
|B(z)|\geq 1/3\quad \text{ for all }\quad z\in \Gamma.
\label{estimacion}
\end{equation}
It is easy to check that the M\"obius transform \ $z\mapsto
(z-\lambda )(z+\lambda)$ \ sends the circle $\Gamma $ onto the \
circle that contains the interval $[-1,(2-\lambda )/(2+\lambda )]$
as a diameter, so that the inequality
\[
\left| \frac{z-\lambda }{z+\lambda }\right| \geq \frac{2-\lambda
}{2+\lambda }=\frac{1-\lambda /2}{1+\lambda /2}
\]
holds for all $z\in \Gamma $, and
\[
|B(z)|\geq \prod_{k=0}^{n}\frac{1-\lambda_{k}/2}{1+\lambda
_{k}/2}.
\]
To estimate the above product, we take into consideration the fact
that for all $x,y\geq 0$, the inequality
\[
\frac{1-x}{1+x}\cdot \frac{1-y}{1+y} \geq \frac{1-(x+y)}{1+x+y}
\]
holds. This leads us to the inequality
\[
\prod_{k=0}^{n}\frac{1-\lambda_{k}/2}{1+\lambda_{k}/2}\geq \frac{1-%
\frac{1}{2}\sum_{k=0}^{n}\lambda
_{k}}{1+\frac{1}{2}\sum_{k=0}^{n}\lambda
_{k}}=\frac{1-1/2}{1+1/2}=1/3\text{,}
\]
which proves (\ref{estimacion}).

It follows from the definition of $T$, that
\[
T^{\prime \prime }(t)=\frac{1}{2\pi \ii}\int_{\Gamma }\frac{%
z^{2}\ee^{-zt}}{B(z)}\dd z
\]
and, using  $\alpha (\theta )=1+\ee^{\ii\theta }$ as a
parametrization of $\Gamma$, we have  (taking into account
Fubini's Theorem) that
\begin{eqnarray*}
\int_0^{\infty}|T^{\prime \prime }(t)| \dd t & \leq & (3/2\pi
)\int_0^{\infty}\int_{0}^{2\pi }2(1+\cos \theta
)\ee^{-(1+\cos \theta )t}\dd \theta \dd t\\
&= &(3/2\pi )\int_{0}^{2\pi }2(1+\cos \theta )\frac{1}{(1+\cos \theta )}%
\dd \theta =6.
\end{eqnarray*}
Now we will compute integrals of the form $\int_{0}^{\infty
}\ee^{-\lambda_{k}t}T^{\prime \prime }(t)\dd t$ in terms of the
scalars $\lambda_{k}$. To do this, we note that
\begin{equation}
\int_{0}^{\infty }\ee^{-\lambda_{k}t}T^{\prime \prime }(t)\dd
t=\frac{1}{2\pi \ii}\int_{\Gamma }\frac{z^{2}}{B(z)(z+\lambda
_{k})}\dd z \label{estimacion2}
\end{equation}
and, taking into consideration the fact that
\[\frac{z^{2}}{B(z)(z+\lambda_{k})%
}\]
 has no poles in the exterior of $\Gamma $, the above integral
depends only on its residue at $\infty $. Now
\[
\frac{1}{z+\lambda_{k}}=\frac{1}{z}\sum_{j=0}^{\infty
}(-1)^{j}(\lambda_{k}/z)^{j}=1/z-\lambda _{k}/z^{2}+\lambda
_{k}^{2}/z^{3}-\cdots
\]
and
\begin{eqnarray*}
\frac{1}{B(z)} &=&\prod_{k=0}^{\infty }\frac{1+\lambda
_{k}/2}{1-\lambda _{k}/2}=1+\frac{2\sum_{k=0}^{\infty }\lambda
_{k}}{z}+\frac{2\left(
\sum_{k=0}^{\infty }\lambda _{k}\right) ^{2}}{z^{2}}+\cdots \\[10pt]
&=&1+\frac{2}{z}+\frac{2}{z^{2}}+\cdots
\end{eqnarray*}
so that
\[
\frac{z^{2}}{B(z)(z+\lambda _{k})}=z+(2-\lambda
_{k})+\frac{\lambda_{k}^{2}-2\lambda_{k}+2}{z}+\cdots\;.
\]
This, in conjunction with (\ref{estimacion2}), leads us to the
formula
\begin{equation}
\int_{0}^{\infty }\ee^{-\lambda_{k}t}T^{\prime \prime }(t)\dd
t=\lambda_{k}^{2}-2\lambda_{k}+2.  \label{estimacion3}
\end{equation}
Now, let $q$ be an exponential polynomial of the form
(\ref{expo}). Then, if we take into consideration
(\ref{estimacion3}), we conclude that
\begin{eqnarray*}
\int_{0}^{\infty }q(t+a)T^{\prime \prime }(t)\dd t
&=&\int_{0}^{\infty }\left(\sum_{k=0}^{n}a_{k}\ee^{-\lambda
_{k}t}\ee^{-\lambda_{k}a}\right)T^{\prime \prime
}(t)\dd t \\[10pt]
&=&\sum_{k=0}^{n}a_{k}\ee^{-\lambda_{k}a}\int_{0}^{\infty
}\ee^{-\lambda
_{k}t}T^{\prime \prime }(t)\dd t \\[10pt]
&=&\sum_{k=0}^{n}a_{k}\ee^{-\lambda_{k}a}(\lambda_{k}^{2}-2\lambda_{k}+2) \\[10pt]
&=&q^{\prime \prime }(a)-2q^{\prime }(a)+2q(a)\text{. }
\end{eqnarray*}
Hence,
\begin{eqnarray*}
\left| q^{\prime \prime }(a)-2q^{\prime }(a)+2q(a)\right|
&=&\left|
\int_{0}^{\infty }q(t+a)T^{\prime \prime }(t)\dd t\right| \\[10pt]
&\leq &\int_{0}^{\infty }\left| q(t+a)T^{\prime \prime }(t)\right|
\dd t\leq \|q\|_{[0,\infty )}\int_{0}^{\infty }\left| T^{\prime
\prime }(t)\right| \dd t
\\[10pt]
&\leq &6\|q\|_{[0,\infty )}\text{ (for all }a\geq 0\text{),}
\end{eqnarray*}
so that
\[
\|q^{\prime \prime }\|_{[0,\infty )}\leq 2\|q^{\prime
}\|_{[0,\infty )}+8\|q\|_{[0,\infty )}.
\]
It is well known that the inequality
\[
\|f^{\prime }\|_{[0,\infty )}^{2}\leq 4\|f\|_{[0,\infty
)}\|f^{^{\prime \prime }}\|_{[0,\infty )}
\]
holds for all functions $f$ $\in C^{(2)}[0,\infty )$ (see \cite
{kolmogorov}), so that
\[
\|q^{\prime }\|_{[0,\infty )}^{2}\leq 4\|q\|_{[0,\infty
)}\|q^{^{\prime \prime }}\|_{[0,\infty )}\leq \|q\|_{[0,\infty
)}(8\|q^{\prime }\|_{[0,\infty )}+16\|q\|_{[0,\infty )})
\]
and
\[
\left( \frac{\|q^{\prime }\|_{[0,\infty )}}{\|q\|_{[0,\infty
)}}\right) ^{2}\leq 8\frac{\|q^{\prime }\|_{[0,\infty
)}}{\|q\|_{[0,\infty )}}+16
\]
which clearly implies that
\[
\frac{\|q^{\prime }\|_{[0,\infty )}}{\|q\|_{[0,\infty )}}\leq 11
\]
for all expressions of the form (\ref{expo}).\endproof

Theorem \ref{newmaninequality} is a nice generalization of the
classical Markov inequality, which states that algebraic
polynomials of degree $\leq n$ (i.e., polynomials of the form
$p(x)=a_0+a_1x+ \cdots+a_nx^n$)  satisfy
\[
\|p'\|_{[-1,1]}\leq n^2\|p\|_{[-1,1]}.
\]
Markov's inequality is indeed related to another classical
inequality, due to Bernstein, which states that algebraic
polynomials of degree $\leq n$ satisfy
\[
|p'(x)|\leq \frac{n}{\sqrt{1-x^2}}\|p\|_{[-1,1]},\qquad\text{ for
all } x\in (-1,1).
\]
It is not by chance that Theorem \ref{newmaninequality} appeared
in the middle of our proof. In fact, the study of certain
classical polynomial inequalities and how they can be extended
(and adapted) from the usual spaces of algebraic polynomials to
M\"{u}ntz spaces, has proven to be a deep tool for studying the
density of these spaces. Concretely, the following generalization
of the classical Bernstein inequality will be of fundamental
importance for the main objectives of this section.


\begin{theorem}[Bounded Bernstein's and \Tchebychev's inequalities]\label{boundedinequalities}
\hskip-1pt Let us assume that $0\leq \lambda _{0}<\lambda _{1}<\cdots $ and $%
\sum_{k=1}^{\infty }1/\lambda _{k}<\infty $. Then for each
$\varepsilon >0$ there are constants $c_{\varepsilon },
c_{\varepsilon }^{\ast }>0$ such that
\[\|p\|_{[0,1]}\leq c_{\varepsilon }\|p\|_{[1-\varepsilon ,1]}\]
 and
\[\|p^{\prime }\|_{[0,1-\varepsilon ]}\leq c_{\varepsilon }^{\ast
}\|p\|_{[1-\varepsilon ,1]}\leq c_{\varepsilon }^{\ast
}\|p\|_{[0,1]}\]
 for all $p\in \Pi (\Lambda).$
\end{theorem}

The proof of this theorem is especially tricky, so that we
postpone it to part 2 of this section. We prefer, at present, to
explain how a clever use of this theorem helps us answer several
questions related to the study of the density of M\"{u}ntz spaces.
In particular, we close this subsection by concluding the proof of
the Full M\"{u}ntz Theorem for the space $C[0,1]$.

\vspace{0.5cm}

\noindent \textbf{Proof of Theorem \ref{muntzc01} when $(\protect\lambda _{k})=(%
\protect\alpha _{k})\cup (\protect\beta _{k})$, \textbf{where} $\protect%
\alpha _{k}\rightarrow 0$ \textbf{and }$\protect\beta
_{k}\rightarrow \infty $.} In this case, the relations
$\sum_{k=1}^{\infty }(\lambda _{k})/(\lambda _{k}^{2}+1)=\infty $
and
\begin{equation}
\sum_{k=1}^{\infty }\alpha _{k}+\sum_{k=1}^{\infty }\frac{1}{\beta _{k}}%
=\infty  \label{cond}
\end{equation}
are equivalent. If (\ref{cond}) holds then we have already proved
that $\Pi(\Lambda )$ is dense in $C[0,1]$. Let us now assume that $%
\sum_{k=1}^{\infty }\alpha _{k}<\infty $ and $\sum_{k=1}^{\infty }%
1/\beta_{k}<\infty $.

Recall that a sequence of functions $(f_k)_{k=0}^n\subseteq C(K)$
is called a Haar system on $K$ if
\[\dim{\,\rm span\,}\{f_k\}_{k=0}^n=n+1\]
 and the only element $f\in
{\rm span\,}\{f_k\}_{k=0}^n$ that vanishes at $n+1$ points is the
zero function.  A special type of Haar systems are \Tchebychev\
systems, which are those given by a sequence of functions
$(f_k)_{k=0}^n\subseteq C(K)$ such that
\[
\mathbf{\det }\left(
\begin{array}{cccc}
f_{0}(x_{0}) & f_{1}(x_{0}) & \cdots & f_{n}(x_{0}) \\
\vdots & \vdots & \ddots & \vdots \\
f_{0}(x_{n}) & f_{1}(x_{n}) & \cdots & f_{n}(x_{n})
\end{array}
\right) >0
\]
holds whenever $x_0<x_1<\cdots < x_n$, $\{x_i\}_{i=0}^n\subseteq
K$.
\begin{proposition}\label{muntzsystem} Let us assume that
$\lambda_0<\lambda_1<\cdots<\lambda_n$. Then
$(x^{\lambda_k})_{k=0}^n$ is a \Tchebychev\ system on
$(0,\infty)$.
\end{proposition}

\Proof. Let
$\Delta=\{(\alpha_0,\alpha_1,\ldots,\alpha_n)\in\mathbb{R}^{n+1}:\exists
i\neq j, \alpha_i=\alpha_j\}$. Then it is not difficult to prove
that
\[
D(\rho_0,\ldots,\rho_n):=\mathbf{\det }\left(
\begin{array}{cccc}
x_{0}^{\rho_0} & x_{0}^{\rho_1} & \cdots & x_{0}^{\rho_n} \\
\vdots & \vdots & \ddots & \vdots \\
x_{n}^{\rho_0} & x_{n}^{\rho_1} & \cdots & x_{n}^{\rho_n}
\end{array}
\right) \neq 0
\]
whenever $(\rho_0,\ldots,\rho_n)\in\mathbb{R}^{n+1}\setminus
\Delta$ and $0<x_0<x_1<\cdots < x_n$. The difficult thing is to
show that this determinant must be positive. Now, we take
$\tau=(\tau_0,\ldots,\tau_n):[0,1]\to\mathbb{R}^{n+1}\setminus
\Delta$ a continuous path such that
$\tau(0)=(\lambda_0,\ldots,\lambda_n)$ and
$\tau(1)=(0,1,\ldots,n)$ (this is possible since
$\lambda_0<\cdots<\lambda_n$). The continuity of $\tau$ implies
that
\[
\text{sign}(D(\tau(0)))=\text{sign}(D(\tau(1)))=+1,
\]
since the last determinant is the well known Vandermonde
determinant (see, e.g., \cite{davis}).
\endproof

Let us assume that $(f_k)_{k=0}^n$ is a \Tchebychev\ system.
Under these conditions, it is possible to prove some interesting
results about uniqueness and characterization of best approximants
from the space ${\rm span\,}\{f_k\}_{k=0}^n$. In particular, the
existence of a unique best approximation to $f_n$ by elements of
${\rm span\,}\{f_k\}_{k=0}^{n-1}$ is guaranteed. If $P_n$ is such
an approximant, then the function $T_{n}:=(f_n-P_n)/\|f_n-P_n\|$
is, by definition, the \Tchebychev\ polynomial associated with the
\Tchebychev\ system $(f_k)_{k=0}^n$. As we shall see, these
polynomials play a main role in our theory. In particular, they
satisfy the following nice interlacing property (for a proof, see
\cite{borweinerdelyi}, page 116):

\begin{theorem}[Zeros of \Tchebychev\ Polynomials]\label{interlacingtheorem}
Let us assume that $\mathcal{T}=(f_0,\ldots,f_{n-1},g)$ and
$\mathcal{S}=(f_0,\ldots,f_{n-1},h)$ are \Tchebychev\ systems on
$[a,b]$ and that $T_n=T_{n,\mathcal{T}}$ and
$S_n=S_{n,\mathcal{S}}$ denote the associated \Tchebychev\
polynomials. If $(f_0,\ldots,f_{n-1},g,h)$ is also a \Tchebychev\
system then the zeros of $T_n$ and $S_n$ interlace (i.e., there
exists exactly one zero of $S_n$ between any two consecutive zeros
of $T_n$).
\end{theorem}

Moreover, the following theorem also holds:

\begin{theorem}[Alternation property of \Tchebychev\ Polynomials]\label{alternation}
Let us assume that $\mathcal{T}=(f_0,\ldots,f_{n-1},f_n)$ is a
\Tchebychev\ system on $[a,b]$ and that $T_n=T_{n,\mathcal{T}}$
denotes the associated \Tchebychev\ polynomial. Then there are
$n+1$ points $x_0<x_1<\cdots<x_{n}$ in $[a,b]$ such that
\[
|T_n(x_i)|=\varepsilon (-1)^i, i=0,1,\ldots,n,
\]
where $\varepsilon\in\{1,-1\}$ is the same for all $i$.
\end{theorem}

Let us use the following notation:

\begin{itemize}
\item  $T_{n,\alpha }$ \ denotes the \Tchebychev\  polynomial
associated to the system $(1,x^{\alpha _{k}})_{k=1}^n$.

\item  $T_{n,\beta }$ denotes the \Tchebychev\ polynomial  associated to the system $%
(1,x^{\beta _{k}})_{k=1}^n$.

\item  $T_{2n,\alpha ,\beta }$ denotes the \Tchebychev\ polynomial
associated to the system
\[(1,x^{\alpha _{1}},x^{\alpha
_{2}},\ldots,x^{\alpha _{n}},x^{\beta _{1}},x^{\beta
_{2}},\cdots,x^{\beta _{n}}).
\]
\end{itemize}

It follows from Newman's inequality
\[
\left\| xT_{n,\alpha }^{\prime }(x)\right\| _{[0,1]}\leq
11M(\{\alpha _{k}\})\left\| T_{n,\alpha }\right\| _{[0,1]}={\rm
const}<\infty \text{ (}n\in \mathbb{N}\text{),}
\]
that for each $\varepsilon >0$ there exists a constant
$k_{1}(\varepsilon )$ which only depends on $\varepsilon $ and
$\{\alpha _{k}\}$ (but does not depend on $n$) such that
$T_{n,\alpha }$ has at most $k_{1}(\varepsilon )$
zeros in $[\varepsilon ,1)$ and at least $n-k_{1}(\varepsilon )$ zeros in $%
[0,\varepsilon )$. (It is not possible to increase the number of
zeros of $\ T_{n,\alpha }$ in $[\varepsilon ,1)$ without
increasing the modulus of the derivative of $T_{n,\alpha }$ at
least in some points of the same interval.) On the other hand, and
due to similar reasons, it follows from the bounded Bernstein's
inequality (Theorem \ref{boundedinequalities}), applied to
$T_{n,\beta }$, that
\[
\|T_{n,\beta }^{\prime }\|_{[0,1-\varepsilon ]}\leq c_{\varepsilon
}^{\ast }\|T_{n,\beta }\|_{[0,1]}=c_{\varepsilon }^{\ast }<\infty.
\]
Hence $T_{n,\beta }$ has at most $k_{2}(\varepsilon )$ zeros in $%
[0,1-\varepsilon )$ and at least $n-k_{2}(\varepsilon )$ zeros in $%
[1-\varepsilon ,1)$.

Now, if we take into account the interlacing properties of the
\Tchebychev's polynomials (Theorem \ref{interlacingtheorem}), and
the fact that the system $(1,x^{\alpha _{k}},x^{\beta
_{k}})_{k=1}^{n}$ is an extension of both systems $(1,x^{\alpha
_{k}})_{k=1}^{n}$ and $(1,x^{\beta _{k}})_{k=1}^{n}$, it follows
that, for $n$ big enough, $T_{2n,\alpha ,\beta }$ has at least
$n-k_1(\varepsilon )-1$ zeros on $[0,\varepsilon]$ and at least
$n-k_2(\varepsilon )-1$ zeros on $[1-\varepsilon,1]$. Hence we
conclude that there exists a certain constant $k=k(\varepsilon )$
(which only depends on the sequence $(\lambda _{k})$) such that
$T_{2n,\alpha ,\beta }$ has at most $k(\varepsilon )$ zeros in the
interval $(\varepsilon ,1-\varepsilon )$.

Set $k=k(1/4)$ and let us take a set of points
\[
1/4<t_{0}<t_{1}<\cdots < t_{k+3}<3/4
\]
and a function $f\in C[0,1]$ such that $f(x)=0$ for all $x\in
\lbrack 0,1/4]\cup \lbrack 3/4,1]$ and $f(t_{i})=(-1)^{i}2$ for
all $0\leq i\leq k+3$. Let us assume that there exists a
polynomial $p\in \Pi (\Lambda ) $ such that $\left\| f-p\right\|
_{[0,1]}<1$. Then $p-T_{2n,\alpha ,\beta } $ has at least $2n+1$
zeros in the interval $(0,1)$ (where we have used that $p$
dominates in $[1/4,3/4]$ and $T_{2n,\alpha ,\beta }$ dominates
outside this interval). This is in contradiction to the fact that $%
p-T_{2n,\alpha ,\beta }\in \Pi _{2n}(\Lambda )$ for all $n$ large
enough (which implies that $p-T_{2n,\alpha ,\beta }$ has at most
$2n$ zeros). This ends the proof whenever $\Lambda $ has no
accumulation points in $(0,\infty )$.\endproof

\noindent \textbf{Proof of Theorem 11 for the case in which
$\Lambda $ has some accumulation point in $(0,\infty )$.} In this
case there
exists an infinite subsequence $(\alpha _{k})\subset \Lambda $ such that $%
\inf \{\alpha _{k}\}>0$ and, in such a case, we already know that
$\Pi ((\alpha _{k}))\subset \Pi (\Lambda )$ is a dense subset of
$C[0,1]$. \endproof

\subsection{Proof of the bounded Bernstein and \Tchebychev\ inequalities}
We devote this subsection to the proof of Theorem
\ref{boundedinequalities}. The proof is long, so we divide it into
several steps:

\subsubsection*{Step 1: Bernstein's and \Tchebychev's exponents
satisfying a jump condition} In this step, we prove Bernstein's
and \Tchebychev's inequalities for sequences of exponents that
satisfy the following jump condition: $\inf_{k\in
\mathbb{N}}(\lambda _{k}-\lambda_{k-1})>0$. In particular, we
start with Bernstein's inequality in this special case:

\begin{theorem}[Bounded Bernstein Inequality for Special Sequences]
\label{bernsteininequality} Let us assume that
$\Lambda=(\lambda _{k})_{k=0}^{\infty }$ is a sequence of
nonnegative real numbers that satisfies the jump condition
$\inf_{k\in \mathbb{N}}(\lambda _{k}-\lambda _{k-1})>0$, and
$\sum_{k=1}^{\infty }1/\lambda _{k}<\infty $, $\lambda _{0}=0,$
$\lambda _{1}\geq 1$. Then for all $\varepsilon \in (0,1)$ there
exists a constant
$c_{\varepsilon }=c(\varepsilon ,\Lambda%
)>0$ such that the inequalities
\[
\|p^{\prime }\|_{[0,1-\varepsilon ]}\leq c_{\varepsilon
}\|p\|_{L^{2}(0,1)}\text{,}\quad \|p^{\prime }\|_{[0,1-\varepsilon
]}\leq c_{\varepsilon }\|p\|_{[0,1]}
\]
hold for all $p\in \Pi (\Lambda)$.
\end{theorem}

\Proof.   It follows from direct computation of the errors
$E(x^{\alpha },\Pi
((\lambda _{k})_{k=0}^{n}))_{2}$ that, for all $m\in \mathbb{N}$ and all $p$ $%
\in \Pi (\Lambda\backslash\{\lambda_m\})$, the inequality
\begin{eqnarray}
\|x^{\lambda _{m}}-p(x)\|_{L^{2}(0,1)} &\geq &\frac{1}{\sqrt{2\lambda _{m}+1}%
}\prod_{k\geq 0\text{; }k\neq m}\left| \frac{\lambda _{m}-\lambda _{k}}{%
\lambda _{m}+\lambda _{k}+1}\right|  \label{+++} \\
&=&\frac{1}{\sqrt{2\lambda _{m}+1}}\prod_{k\geq 0\text{; }k\neq
m}\left| \frac{(\lambda _{k}+1/2)-(\lambda _{m}+1/2)}{(\lambda
_{k}+1/2)+(\lambda _{m}+1/2)}\right|  \nonumber
\end{eqnarray}
holds. Hence it is of interest to study products of the form:
\[
\prod_{k\geq 0\text{; }k\neq m}\left| \frac{\alpha_{k}+\alpha
_{m}}{\alpha _{k}-\alpha _{m}}\right|
\]
for sequences $(\alpha _{k})_{k=0}^{\infty }$ such that
$\inf_{k\in \mathbb{N}}(\alpha _{k}-\alpha_{k-1})>0$, and
$\sum_{k=0}^{\infty }1/\alpha _{k}<\infty $. (Note that we have,
for ease of exposition, reversed the quotients.)

We decompose:
\begin{eqnarray*}
\prod_{k\geq 0\text{; }k\neq m}\left| \frac{\alpha _{k}+\alpha _{m}}{%
\alpha _{k}-\alpha _{m}}\right|
&=&\prod_{k\geq 0\text{; }\alpha _{k}<\alpha _{m}}\left| 1+\frac{2\alpha _{m}%
}{\alpha _{k}-\alpha _{m}}\right|  \times\\
&\times&\prod_{k\geq 0\text{; }\alpha _{m}<\alpha
_{k}<2\alpha _{m}}\left| 1+\frac{2\alpha _{m}}{\alpha _{k}-\alpha _{m}}%
\right| \prod_{k\geq 0\text{; }\alpha _{k}\geq 2\alpha _{m}}\left| 1+\frac{%
2\alpha _{m}}{\alpha_{k}-\alpha _{m}}\right|.
\end{eqnarray*}
Clearly, for all $k$ such that $\alpha _{k}\geq 2\alpha _{m}$ we have that $%
\alpha _{k}/2\geq \alpha _{m}$ so that $\alpha_{k}-\alpha _{m}\geq
\alpha_{k}-\alpha_{k}/2=\alpha _{k}/2$. This means that
\begin{eqnarray*}
\prod_{k\geq 0\text{; }\alpha _{k}\geq 2\alpha _{m}}\left|
1+\frac{2\alpha
_{m}}{\alpha _{k}-\alpha_{m}}\right| &\leq &\exp\! \left( \sum_{k\geq 0\text{%
; }\alpha_{k}\geq 2\alpha_{m}}\left| \frac{2\alpha _{m}}{\alpha
_{k}-\alpha _{m}}\right| \right) \\
&\leq &\exp\! \left( 4\alpha _{m}\sum_{k\geq 0\text{; }\alpha
_{k}\geq 2\alpha _{m}}\frac{1}{\alpha_{k}}\right),
\end{eqnarray*}
which implies that there exists a constant $\xi _{m}>0$ such that
$\xi _{m}\leq 4\sum_{k\geq 0\text{; }\alpha _{k}\geq 2\alpha
_{m}}1/\alpha_{k}$ and
\[
\prod_{k\geq 0\text{; }\alpha_{k}\geq 2\alpha _{m}}\left|
1+\frac{2\alpha _{m}}{\alpha_{k}-\alpha _{m}}\right| =\exp
(\alpha_{m}\xi _{m}).
\]
The products
\[
\prod_{k\geq 0\text{; }\alpha_{k}<\alpha _{m}}\left| 1+\frac{2\alpha _{m}}{%
\alpha_{k}-\alpha _{m}}\right|
\]
are bounded. Moreover, we can use that in between $\alpha_m$ and
$2\alpha_m$ there are $o(\alpha_m)$ terms (since
$\sum_{k}\frac{1}{\alpha_k}<\infty$) to estimate the product
$\prod_{k\geq 0\text{; }\alpha _{m}<\alpha_{k}<2\alpha _{m}}\left| 1+\frac{%
2\alpha_{m}}{\alpha_{k}-\alpha _{m}}\right|$ and prove that there
are constants $\gamma _{m}$ such that
\[
\prod_{k\geq 0\text{; }k\neq m}\left| \frac{\alpha_{k}+\alpha
_{m}}{\alpha _{k}-\alpha _{m}}\right| \le\exp (\alpha _{m}\gamma
_{m}), \qquad \lim_{m\rightarrow \infty }\gamma _{m}=0.
\]
 Hence
\[
\prod_{k\geq 0\text{; }k\neq m}\left| \frac{\alpha _{k}-\alpha
_{m}}{\alpha _{k}+\alpha _{m}}\right| \ge \exp (-\alpha _{m}\gamma
_{m}), \qquad \lim_{m\rightarrow \infty }\gamma _{m}=0
\]
and, taking into consideration the formula (\ref{+++}), we obtain
that
\[
\|x^{\lambda _{m}}-p(x)\|_{L^{2}(0,1)}\geq \exp (-\gamma _{m}\lambda _{m})%
\]
where $\lim_{m\rightarrow \infty }\gamma _{m}=0$ and
$p\in \Pi((\Lambda\backslash\{\lambda_m\})$. This clearly implies
that for every polynomial $p=\sum a_{k}x^{\lambda _{k}}\in \Pi
(\Lambda)$, the inequality
\[
\|p\|_{L^{2}(0,1)}=\|a_{k}x^{\lambda _{k}}-(a_{k}x^{\lambda
_{k}}-p(x))\|_{L^{2}(0,1)}\geq |a_{k}|\exp (-\gamma _{k}\lambda
_{k})
\]
holds. Hence
\begin{equation}
|a_{k}|\leq \exp (\gamma _{k})^{\lambda
_{k}}\|p\|_{L^{2}(0,1)}\leq c_{\varepsilon }(1+\varepsilon
)^{\lambda _{k}}\|p\|_{L^{2}(0,1)} \label{des1}
\end{equation}
for a certain constant $c_{\varepsilon }$, since only a finite
number of values $\gamma _{k}$ satisfy $\exp (\gamma
_{k})>1+\varepsilon $ (the constant $c_{\varepsilon }$ only
depends on the behaviour of the other values $\gamma _{k}$).
Another proof -- indeed the original one -- of this inequality was
given by Clarkson and Erdös \cite{clarkson} in 1943.

Taking into consideration that $\lambda _{1}\geq 1$ and
$c=\inf_{k\in \Bbb{N}} (\lambda _{k}-\lambda _{k-1})>0$, we have
that there exists a strictly
increasing sequence of natural numbers $m_{j}$, $j\geq 0$, such that $%
\{\floor{\lambda _{k}}\}_{k=0}^{\infty }=\{m_{j}\}_{j=0}^{\infty }$, where $%
\floor{\lambda _{k}}$ denotes the integer part of $\lambda _{k}$
for each $k\in \mathbb{N}$. Furthermore,
\[
M:=M(\Lambda):=\max_{j\geq 0}\#\{k:\floor{\lambda
_{k}}=m_j\}<\infty,
\]
so that
\begin{eqnarray*}
\sum_{k=0}^{n}\lambda _{k}\alpha ^{\lambda _{k}-1} &\leq
&\sum_{k=0}^{n}(\floor{\lambda _{k}}+1)\alpha ^{\floor{\lambda _{k}}-1} \\
&=&\sum_{k=0}^{n}\floor{\lambda _{k}}\alpha ^{\floor{\lambda
_{k}}-1}+\sum_{k=0}^{n}\alpha ^{\floor{\lambda _{k}}-1} \\
&\leq &M\left( \sum_{k=0}^{n}m_{k}\alpha
^{m_{k}-1}+\sum_{k=0}^{n}\alpha
^{m_{k}-1}\right) \\[5pt]
&\leq &M\left( \sum_{k=0}^{\infty }m_{k}\alpha
^{m_{k}-1}+\sum_{k=0}^{\infty
}\alpha ^{m_{k}-1}\right) \\[5pt]
&\leq &C(\alpha,M):=M\left( \sum_{t=1}^{\infty }t\alpha
^{t-1}+\sum_{t=1}^{\infty }\alpha ^{t-1}\right) <\infty
\end{eqnarray*}
for all $\alpha \in (0,1)$.

We can use the above inequality to estimate the norm
$\|p^{\prime}\|_{[0,1-\varepsilon ]}$ as follows:
\begin{eqnarray*}
\|p^{\prime }\|_{[0,1-\varepsilon ]} &\leq &\sum_{k=0}^{n}\lambda
_{k}|a_{k}|(1-\varepsilon )^{\lambda _{k}-1} \\
&\leq &c_{\varepsilon }(1+\varepsilon )\sum_{k=0}^{n}\lambda
_{k}[(1+\varepsilon )(1-\varepsilon )]^{\lambda
_{k}-1}\|p\|_{L^{2}(0,1)}
\\
&=&c_{\varepsilon }(1+\varepsilon )\sum_{k=0}^{n}\lambda
_{k}[(1-\varepsilon
^{2})]^{\lambda _{k}-1}\|p\|_{L^{2}(0,1)} \\
&\leq &c_{\varepsilon }(1+\varepsilon )C(1-\varepsilon
^{2},M)\|p\|_{L^{2}(0,1)} \\[5pt]
&\leq &c(\varepsilon ,\Lambda )\|p\|_{L^{2}(0,1)}\leq
c(\varepsilon ,\Lambda )\|p\|_{C[0,1]}\text{,}
\end{eqnarray*}
which is what we wanted to prove.\endproof

The next theorem, proved by L. Schwartz \cite{schwartz}, is an
important consequence of the inequality (\ref{des1}).

\begin{theorem}[Closure of Nondense M\"{u}ntz Spaces for Special Sequences]
\label{closurespecial} Let us assume that $%
\Lambda=(\lambda _{k})_{k=0}^{\infty }$ is a sequence of
nonnegative real numbers
such that $\inf_{k\in \mathbb{N}}(\lambda _{k}-\lambda _{k-1})>0$ \ and $%
\sum_{k=1}^{\infty }1/\lambda _{k}<\infty $, $\lambda _{0}=0,$
$\lambda _{1}\geq 1$. Then the functions that belong to the
closure of $\Pi(\Lambda)$ can be analytically extended to $\mathbb{D}%
\setminus [-1,0]$. If $\lambda_{k}$ is an integer for all $k$%
, then the functions of the closure of $\Pi (\Lambda)$ can be
analytically extended to the unit disc. Finally, if $\Lambda$ is
lacunary $($i.e., $\inf\{\lambda _{k}/\lambda
_{k-1}\}_{k=2}^{\infty }> 1)$  then the closure of $\Pi (\Lambda)$
is precisely the set
\[
\left\{ f\in C[0,1]:f(x)=\sum_{k=0}^{\infty }a_{k}x^{\lambda _{k}},%
\text{ }x\in \lbrack 0,1]\right\} .
\]
\end{theorem}

\Proof.  Let us assume that $\lim_{n\rightarrow \infty }\|f-q_{n}\|_{C%
[0,1]}=0$, where $q_{n}(x):=\sum_{k=0}^{k_{n}}a_{n,k}x^{\lambda
_{k}}$. Then the sequence of polynomials $(q_{n})_{n=0}^{\infty }$
is a Cauchy sequence in $C[0,1]$. It follows that for each $\delta
>0$ and all $n\in \mathbb{N}$,
\[
|a_{n,k}-a_{m,k}|\leq c_{\delta }(1+\delta )^{\lambda _{k}}\|q_{n}-q_{m}\|_{%
C[0,1]}\rightarrow 0\text{ ( }n,m\rightarrow \infty \text{). }
\]
This means that there are numbers $a_{k}\in \Bbb{R}$ such that $%
\lim_{n\rightarrow \infty }a_{n,k}=a_{k}$ ($k\in \mathbb{N}$). Let $%
h(x):=\sum_{k=0}^{\infty }a_{k}x^{\lambda _{k}}$. Then for all
$\delta >0$ we can write
\[
|a_{k}|=\lim_{n\rightarrow \infty }|a_{n,k}|\leq
\lim_{n\rightarrow \infty }c_{\delta }(1+\delta )^{\lambda
_{k}}\|q_{n}\|_{C[0,1]}=c_{\delta }(1+\delta )^{\lambda
_{k}}\|f\|_{C[0,1]}.
\]

It follows that the series $h(x)=\sum_{k=0}^{\infty
}a_{k}x^{\lambda _{k}}$ is absolutely convergent for all $x<1$. To
prove this claim we take into account that $\lambda _{k}/k\geq c$
for all $k$, so that the inequality
\[
|a_{k}x^{\lambda _{k}}|^{1/k}\leq (c_{\delta }\|f\|_{C%
[0,1]})^{1/k}((1+\delta )x)^{\lambda _{k}/k}\leq (c_{\delta }\|f\|_{{C%
}[0,1]})^{1/k}((1+\delta )x)^{c}< 1
\]
holds for $k$ sufficiently large and $\delta $ such that
$(1+\delta )x<1$. Now, it is clear that $h$ coincides with the
function $f$. Consider the branch of logarithm that is defined on
the complex plane cut along $(-\infty,0]$ and that is positive for
values $>1$. For any $\lambda$ this defines a branch of
$z^\lambda=\exp(\lambda\log z)$. Now
 if $z\in \mathbb{D}\setminus \lbrack -1,0]$ then
\[
\sum_{k=0}^{\infty }\left|a_{k}z^{\lambda _{k}}\right| \leq
f(|z|)<\infty.
\]
This proves that $f(z)=\sum_{k=0}^{\infty }a_{k}z^{\lambda _{k}}$
is analytic on $\mathbb{D}\setminus \lbrack -1,0]$. If $\lambda
_{k}$ is an integer for all $k$, then it is clear that the above
arguments are also valid for all $z\in (-1,0)$ so that $f$ is
analytic in the unit disc. If the sequence is lacunary then we can
use a theorem by Hardy and Littlewood \cite{hardy} that claims
that if the power series $\sum_k a_kx^\lambda_k$, with radius of
convergence $1$, is lacunary and $\lim_{x\to
1^-}\sum_ka_kx^{\lambda_k}=a$, then $\sum_ka_k=a$, to conclude
that the series $\sum_{k=0}^{\infty }a_{k}z^{\lambda _{k}}$ also
converges for $z=1$. In the other cases there are counterexamples,
i.e., there are series of the form $\sum_{k=0}^{\infty
}a_{k}z^{\lambda _{k}}$ that belong to the closure of $\Pi
(\Lambda)$ in $C [0,1]$ and they do not converge for $z=1$ (see
\cite{clarkson}).
\endproof

We prove, for the special case we  are considering in this step,
\Tchebychev's inequality which claims that  the norms of the
elements of (nondense)  M\"{u}ntz spaces essentially depend on the
behaviour of the elements near $x=1$.

\begin{corollary}[Bounded \Tchebychev\ Inequality for Special Sequences]
\label{boundedchebychevineq}\hskip-5pt Under the
hypotheses of Theorem $\ref{bernsteininequality}$, for each
$\varepsilon \in (0,1)$ there exists a constant $c_{\varepsilon
}=c(\varepsilon ,\Lambda)$ such that $\|p\|_{C[0,1]}\leq
c_{\varepsilon }\|p\|_{C[1-\varepsilon ,1]}$ for all $p\in \Pi
(\Lambda)$.
\end{corollary}

\Proof.  Making (if  necessary) the change of variable
$y=x^{1/\lambda _{1}}$ we may assume, without loss of generality,
that $\lambda _{1}=1$. Let us now assume that there
exists a sequence of polynomials $(p_{n})\subset
\Pi (\Lambda)$ such that $\lim_{n\rightarrow \infty }\|p_{n}\|_{%
C[0,1]}=\infty $ but $\|p_{n}\|_{C[1-\varepsilon ,1]}=1$ for all $n$%
. Then $q_{n}:=p_{n}/\|p_{n}\|_{C[0,1]}$ satisfies $\|q_{n}\|_{%
C[0,1]}=1$ for all $n\in \mathbb{N}$, and $\lim_{n\rightarrow
\infty }\|q_{n}\|_{C[1-\varepsilon ,1]}=0$. It follows from the
bounded Bernstein inequality that for each $\delta \in (0,1)$
there exists a constant $c_{\delta }$ such that $\|q_{n}^{\prime
}\|_{[0,1-\delta ]}\leq c_{\delta }$
for all $n$. We may use the Arzel\`{a}-Ascoli theorem in the interval $%
[0,1-\varepsilon /2]$ to \ obtain from $(q_{n})$ a subsequence
that converges uniformly to a certain $f\in C[0,1-\varepsilon
/2]$. Using the same arguments as in the proof of Theorem
\ref{closurespecial}, we get more information: $f$ must be
analytic on $(0,1-\varepsilon /2)$. But $\lim_{n\rightarrow \infty
}\|q_{n}\|_{C[1-\varepsilon ,1]}=0$ implies that
$f|_{(1-\varepsilon,1-\varepsilon  /2)}=0$, which clearly implies
that $f$ is the null function (just apply the well known Identity
Principle of complex analysis). The fact that $f=0$ and
$\|q_{n}\|_{C[0,1]}=1$ for all $n$ are simultaneously impossible.
\endproof

\subsubsection*{Step 2. Comparison results }

The main goal of this step is to introduce a few results that
will be useful for the proof, the next step, of Bernstein's and
\Tchebychev's inequalities for general sequences of exponents
$(\lambda_k)_{k=0}^{\infty}$. These results are expressed in terms
of the \Tchebychev\ polynomials associated with the M\"{u}ntz
system $(x^{\lambda_k})_{k=0}^{\infty}$.

Let us proceed by stages. We first introduce some notation. We say
that $(f_{0},\ldots,f_{n})$ is a Descartes system on $[a,b]$ if
\[
\mathbf{\det }\left(
\begin{array}{cccc}
f_{i_{0}}(x_{0}) & f_{i_{1}}(x_{0}) & \cdots & f_{i_{m}}(x_{0}) \\
\vdots & \vdots & \ddots & \vdots \\
f_{i_{0}}(x_{m}) & f_{i_{1}}(x_{m}) & \cdots & f_{i_{m}}(x_{m})
\end{array}
\right) >0
\]
holds whenever $0\leq i_{0}<i_{1}<\cdots <i_{m}\leq n$ and $a\leq
x_{0}<x_{1}<\cdots <x_{m}\leq b$.  We say that
$(f_{0},\ldots,f_{n})$  is a Markov system in $C[a,b]$ if for all
$f\in C[a,b] $ and all $k\leq n$ there exists a unique best
approximation to $f$ by elements of $M_{k}:={\rm
span\,}\{f_{i}\}_{i=0}^{k}$. The
$n$th \Tchebychev\ polynomial associated with the Markov system  $(f_{i})%
_{i=0}^{n}$ is given by $%
T_{n}:=(f_{n}-p_{n-1})/\|f_{n}-p_{n-1}\|_{C[a,b]}$, where $p_{n-1}$ is the unique best approximation to $%
f_{n}$ by elements of $M_{n-1}.$ Sometimes, by a misuse of
notation, we also say that the \Tchebychev\  polynomial $T_n$ is
associated with $M_{n}$.

Now, one of the main properties of M\"{u}ntz spaces is that
\[
(x^{\lambda _{0}},x^{\lambda _{1}},\ldots,x^{\lambda _{n}}),
\]
where $0=\lambda _{0}<\lambda _{1}<\cdots <\lambda _{n}$, is a
Descartes system on each interval $[a,b]\subset [0,\infty)$ (see
Proposition \ref{muntzsystem}). The following technical lemma,
whose proof is quite involved, is a nice refinement of the
classical Descartes rule of signs and was proved by Pinkus and
(independently) by P. W. Smith. It is the key for the proof of the
comparison results we will need (see \cite[p.~103]
{borweinerdelyi}, or \cite{smith} for a proof).

\begin{lemma}[Pinkus-Smith] Let us assume that $(f_{0},\ldots,f_{n})$
is \ a Descartes system on $[a,b]$, and let
\[
p=f_{k}+\sum_{i=1}^{r}a_{i}f_{k_{i}}\text{; }q=f_{k}+%
\sum_{i=1}^{r}b_{i}f_{t_{i}};\quad\text{ with }\ a_{i},b_{i}\in
\Bbb{R}
\]
be chosen such that $0\leq t_{i}\leq k_{i}<k$ for all $i\in
\{1,\ldots,m\}$ and $k<t_{i}\leq k_{i}\leq n$ for all $i\in
\{m+1,\ldots,r\}$, with strict inequality for at least one of the
indices $i\in \{1,\ldots,r\}$.

\noindent If $p(x_{i})=q(x_{i})=0$ for the distinct points
$x_{i}\in \lbrack a,b]$, $i=1,\ldots,r$, then
\[
|p(x)|\leq |q(x)|, \qquad x\in \lbrack a,b].
\]
Furthermore, the inequality is strict for all $x\in \lbrack
a,b]\setminus \{x_{i}\}_{i=1}^{r}$.
\end{lemma}

We use this result with the M\"{u}ntz spaces ${M}_{n}(\Lambda
)=\Pi ((\lambda _{k})_{k=0}^{n})$ and ${M}_{n}(\Gamma )=\Pi
((\gamma _{k})_{k=0}^{n})$, where we assume that $0=\lambda
_{0}<\lambda
_{1}<\cdots <\lambda _{n}$, $0=\gamma _{0}<\gamma _{1}<\cdots <\gamma _{n}$, and $%
\lambda _{k}\geq \gamma _{k}$ for all $k$. With this idea in mind,
we take $s\in (0,1)$ and denote by
$T_{n,\lambda }$ and $T_{n,\gamma }$ the \Tchebychev\  polynomials associated with $%
{M}_{n}(\Lambda )$ and ${M}_{n}(\Gamma )$, respectively, on the
interval $[1-s,1]$.

\begin{lemma}
With the hypotheses and notation just introduced, the following
claims hold:
\begin{itemize}
\item[$(a)$]  Let $y\in \lbrack 0,1-s)$. Then the maximum values
of the expressions
\[
\max_{0\neq p\in {M}_{n}(\Lambda )}\frac{|p(y)|}{\|p\|_{[1-s,1]}}%
\qquad\text{ and }\qquad \max_{0\neq p\in {M}_{n}(\Lambda )}\frac{|p^{\prime }(y)|}{%
\|p\|_{[1-s,1]}}
\]
are both attained by $p=T_{n,\lambda }$. (In the second case we assume that $%
\lambda _{1}\geq 1$ whenever $y=0$.)

\item[$(b)$]  $|T_{n,\lambda }(0)|\leq |T_{n,\gamma }(0)|$. Furthermore, if $%
\lambda _{1}=\gamma _{1}=1$ then also $|T_{n,\lambda }^{\prime
}(0)|\leq |T_{n,\gamma }^{\prime }(0)|$.
\end{itemize}
\end{lemma}

\Proof.  We propose the proof of $(a)$ as an exercise. Let us
prove $(b)$. Let $p\in M_{n}(\Gamma )$ be such that it
interpolates $T_{n,\lambda }$ at its zeros (which are all of them
simple
zeros), and in $(0,1)$. It follows from the Pinkus-Smith Lemma that $|p(x)|\leq $ $%
|T_{n,\lambda }(x)|$ for all $x\in \lbrack 0,1]$. In particular, $%
\|p\|_{[1-s,1]}\leq \|T_{n,\lambda }\|_{[1-s,1]}=1$ and, taking
into account part (a) of this lemma, we get
\[
|T_{n,\lambda }(0)|=|p(0)|\leq \frac{|p(0)|}{\|p\|_{[1-s,1]}}\leq \frac{%
|T_{n,\gamma }(0)|}{\|T_{n,\gamma }\|_{[1-s,1]}}=|T_{n,\gamma
}(0)|,
\]
which proves the first part of (b). To prove the second claim, the
argument is similar. We take $0\neq p\in M_{n}(\Gamma )$ such that
it interpolates $T_{n,\lambda }$ at its zeros in $[1-s,1]$ (there
are $n$ zeros), and we normalize by imposing the additional
condition $p^{\prime }(0)=T_{n,\lambda }^{\prime }(0)$. (Note that
$p^{\prime }(0)\neq 0$, since otherwise we would have that $p\in
{\rm span\,}\{x^{\gamma _{k}}:k =
0,2,3,\ldots,n\}$ has $n$ zeros in $[1-s,1]$, which is impossible since $%
(x^{\gamma _{k}}:k = 0,2,3,\ldots,n)$ is a Descartes system.) Then $%
|p(x)|\leq $ $|T_{n,\lambda }(x)|$ for all $x\in \lbrack 0,1]$. Hence $%
\|p\|_{[1-s,1]}\leq \|T_{n,\lambda }\|_{[1-s,1]}=1$ and it follows
again from part (a) of this lemma that
\[
|T_{n,\lambda }^{\prime }(0)|=|p^{\prime }(0)|\leq \frac{|p^{\prime }(0)|}{%
\|p\|_{[1-s,1]}}\leq \frac{|T_{n,\gamma }^{\prime
}(0)|}{\|T_{n,\gamma }\|_{[1-s,1]}}=|T_{n,\gamma }^{\prime
}(0)|\text{,}
\]
which proves the second part of (b).\endproof

\begin{lemma}
$|T_{n,\lambda }(x)|$ and $|T_{n,\gamma }(x)|$ are monotone
decreasing functions on the interval $[0,1-s]$. Furthermore, if
$\lambda _{1}=\gamma _{1}=1$, then also $|T_{n,\lambda }^{\prime
}(x)|$ and $|T_{n,\gamma }^{\prime }(x)|$ are monotone decreasing
on the interval $[0,1-s]$.
\end{lemma}

\Proof.  Let us assume that $|T_{n,\lambda }(x)|$ is not monotone decreasing on $%
[0,1-s]$. Then $T_{n,\lambda }^{\prime }(x)\in {\rm
span\,}\{x^{\lambda _{k}-1}: k\in \{1,2,3,\ldots,n\}\}$ has at
least $n$ zeros in $(0,1)$, which is impossible. The second claim
can be proved by similar arguments.\endproof

\begin{theorem}[Comparison Theorem]\label{comparisontheorem} The inequality
\[
\max_{0\neq p\in M_{n}(\Lambda )}\frac{\|p\|_{[0,1]}}{%
\|p\|_{[1-s,1]}}\leq \max_{0\neq p\in M_{n}(\Gamma )}\frac{%
\|p\|_{[0,1]}}{\|p\|_{[1-s,1]}}
\]
holds. Furthermore, if $\lambda _{1}=\gamma _{1}=1$ then
\[
\max_{0\neq p\in M_{n}(\Lambda )}\frac{\|p^{\prime }\|_{[0,1-s]}}{%
\|p\|_{[1-s,1]}}\leq \max_{0\neq p\in M_{n}(\Gamma )}\frac{%
\|p^{\prime }\|_{[0,1-s]}}{\|p\|_{[1-s,1]}}.
\]
\end{theorem}

\Proof.   Let $y\in \lbrack 0,1-s)$. Then
\begin{eqnarray*}
\max_{0\neq p\in M_{n}(\Lambda )}\frac{|p(y)|}{\|p\|_{[1-s,1]}} &=&%
\frac{|T_{n,\lambda }(y)|}{\|T_{n,\lambda
}\|_{[1-s,1]}}=|T_{n,\lambda
}(y)|\leq |T_{n,\gamma }(0)| \\
&=&\frac{|T_{n,\gamma }(0)|}{\|T_{n,\lambda }\|_{[1-s,1]}}\leq
\max_{0\neq
p\in M_{n}(\Gamma )}\frac{\|p\|_{[0,1-s]}}{\|p\|_{[1-s,1]}} \\
&\leq &\max_{0\neq p\in M_{n}(\Gamma )}\frac{\|p\|_{[0,1]}}{%
\|p\|_{[1-s,1]}}.
\end{eqnarray*}
On the other hand, if $y\in \lbrack 1-s,1]$, then
\[
\max_{0\neq p\in M_{n}(\Lambda
)}\frac{|p(y)|}{\|p\|_{[1-s,1]}}\leq
1\leq \max_{0\neq p\in M_{n}(\Gamma )}\frac{\|p\|_{[0,1]}}{%
\|p\|_{[1-s,1]}}.
\]
Hence
\[
\max_{0\neq p\in M_{n}(\Lambda )}\frac{\|p\|_{[0,1]}}{%
\|p\|_{[1-s,1]}}\leq \max_{0\neq p\in M_{n}(\Gamma )}\frac{%
\|p\|_{[0,1]}}{\|p\|_{[1-s,1]}}\text{,}
\]
which is what we wanted to prove. By analogous arguments, we have
that
\begin{eqnarray*}
\max_{0\neq p\in M_{n}(\Lambda )}\frac{|p^{\prime }(y)|}{%
\|p\|_{[1-s,1]}} &=&\frac{|T_{n,\lambda }^{\prime
}(y)|}{\|T_{n,\lambda }\|_{[1-s,1]}}=|T_{n,\lambda }^{\prime
}(y)|\leq |T_{n,\lambda }^{\prime
}(0)|\leq |T_{n,\gamma }^{\prime }(0)| \\[10pt]
&=&\frac{|T_{n,\gamma }^{\prime }(0)|}{\|T_{n,\gamma
}\|_{[1-s,1]}}\leq
\max_{0\neq p\in M_{n}(\Gamma )}\frac{\|p^{\prime }\|_{[0,1-s]}}{%
\|p\|_{[1-s,1]}}\text{,}
\end{eqnarray*}
which is the second claim of the theorem.\endproof

\begin{remark} \label{remarkoncomparisontheorem}
{\sl \ It is possible $($with similar arguments$)$ to extend  Theorem
$\ref{comparisontheorem}$ to include M\"{u}ntz polynomials with
arbitrary real exponents $($i.e., we can also consider negative
powers of $x)$.
}
\end{remark}

\subsubsection*{Step 3. The General Bernstein and \Tchebychev\
Inequalities }
 We now complete the proof of Theorem
\ref{boundedinequalities}.

We know that $\lim_{k\rightarrow \infty }\lambda _{k}/k=\infty $ ,
since $\sum_{k=1}^{\infty }1/\lambda _{k}$
converges and $(\lambda_k)$ is monotone. Let $m\in \mathbb{N}$ be
such that $\lambda _{k}>2k$ for all $k\geq m$, and let us take
$\Gamma :=(\gamma _{k})_{k=0}^{\infty }$ defined by:
\[
\gamma _{k}:=\left\{
\begin{array}{lll}
\min \{\lambda _{k},k\}, &  & \text{if }k\in \{0,1,\ldots,m\}\text{,} \\[10pt]
\frac{1}{2}\lambda _{k}+k, &  & \text{if }k>m\text{.}
\end{array}
\right.
\]
Then $\sum_{k=1}^{\infty }1/\gamma _{k}<\infty $, $0\leq \gamma
_{0}<\gamma _{1}<\cdots$, and
\[
\gamma _{k}-\gamma _{k-1}=\left\{
\begin{array}{lll}
\min \{\lambda _{k},k\}-\min \{\lambda _{k-1},k-1\}, &  & \text{if
}k\in
\{0,1,\ldots,m\} \\[10pt]
\frac{1}{2}\lambda _{m+1}+m+1-\min \{\lambda _{m},m\}, &  &
\text{if }i=m+1
\\[10pt]
\frac{1}{2}(\lambda _{k}-\lambda _{k-1})+1, &  & \text{if }k>m+1
\end{array}
\right.
\]
satisfies $ \gamma _{k}-\gamma _{k-1} \geq 1$ for all $k\in
\mathbb{N}$. Furthermore $\gamma _{k}\leq \lambda _{k}$ for all
$k$. This implies (using Theorems \ref{bernsteininequality} and
\ref{comparisontheorem}, Corollary
 \ref{boundedchebychevineq} and Remark \ref{remarkoncomparisontheorem}) that the inequalities
\[
\max_{0\neq p\in M_{n}(\Lambda )}\frac{\|p\|_{[0,1]}}{%
\|p\|_{[1-s,1]}}\leq \max_{0\neq p\in M_{n}(\Gamma )}\frac{%
\|p\|_{[0,1]}}{\|p\|_{[1-s,1]}}=c_{\varepsilon }<\infty
\]
and
\[
\max_{0\neq p\in M_{n}(\Lambda )}\frac{\|p^{\prime
}\|_{[0,1-\varepsilon ]}}{\|p\|_{[1-s,1]}}\leq \max_{0\neq p\in M%
_{n}(\Gamma )}\frac{\|p^{\prime }\|_{[0,1-\varepsilon ]}}{\|p\|_{[1-s,1]}}%
=c_{\varepsilon }^{\ast }<\infty \text{,}
\]
both hold. \endproof

\begin{corollary}\label{Muntzandrelativecompactness}
Let $(p_{n})_{n=0}^{\infty }$ be a sequence of polynomials in $\
\Pi (\Lambda)$ uniformly bounded in $[0,1]$, and let us assume
that $0\leq \lambda _{0}<\lambda _{1}<\cdots $ and
$\sum_{k=1}^{\infty
}1/\lambda _{k}<\infty $. Then for each $a\in (0,1)$, the sequence $%
(p_{n})_{n=0}^{\infty }$ is a relatively compact subset of $C%
[0,a] $.
\end{corollary}

\Proof.  Set $\varepsilon =1-a$. Then $\|p_{n}^{\prime
}\|_{[0,a]}\leq c_{\varepsilon }^{\ast }\|p_{n}\|_{[a,1]}\leq
c_{\varepsilon }^{\ast }M$ for
all $n$, where $M=\sup \|p_{n}\|_{[0,1]}<\infty $. This implies that $%
\{p_{n}\}_{n=0}^{\infty }$ is equicontinuous in $C[0,a].$ The
corollary follows from the well known Arzel\`{a}-Ascoli
Theorem.\endproof

\subsection{Description of the closure of nondense M\"{u}ntz spaces: the
C[0,1] case}

If $\Pi (\Lambda )$ is not a dense subspace of $C[0,1],$ it is
natural to ask what is its topological closure. Since the density
of $\Pi(\Lambda)$ depends on the convergence character of a
certain series associated with the sequence of exponents
$\Lambda$, it is clear that given two nondense M\"{u}ntz spaces
$\Pi(\Lambda_1)$, $\Pi(\Lambda_2)$, their sum
$\Pi(\Lambda_1)+\Pi(\Lambda_2)=\Pi(\Lambda_1\cup\Lambda_2)$ is
also nondense in $C[0,1]$. This means that the closures of
nondense M\"{u}ntz spaces $\Pi(\Lambda)$ are, in a certain sense,
of small dimension, when viewed as subspaces of $C[0,1]$.

This observation was first made in a famous paper by Clarkson and
Erd\H{o}s \cite{clarkson} published in 1943 in the Duke
Math.~Journal. They proved, for the case of integer exponents
$\Lambda=(n_k)_{k=0}^{\infty}\subset \mathbb{N}$, that
$\sum_{k=0}^{\infty}1/n_k<\infty$ implies that the elements in the
closure of $\Pi(\Lambda)$ are analytic functions defined inside
the unit circle and that their Maclaurin series involves only the
powers $x^{n_k}$ and may diverge at the point $z=1$. Moreover, if
the sequence of exponents is lacunary (which means that
$\inf_{k\geq 0}n_{k+1}/n_k=c>1$), this series converges for $z=1$.
Finally, they used this result to prove, in the particular case
where the exponents are nonnegative integers and for intervals
away from the origin (i.e., intervals $[a,b]$ with $0\not\in
[a,b]$), the natural extension of M\"{u}ntz' theorem (i.e., they
proved that $\sum_{k=0}^{\infty}1/n_k=\infty$ is the necessary and
sufficient condition for density of the M\"{u}ntz space
$\Pi((n_k)_{k=0}^{\infty})$ independently of the appearance or not
of the zero power  in the exponents sequence
$(n_k)_{k=0}^{\infty}$).

This same question was tackled by L.~Schwartz \cite{schwartz} for
certain strictly increasing sequences of exponents (he assumed
$\inf_{k\in\mathbb{Z}}(\lambda_k-\lambda_{k-1})>0$ and proved
Theorem \ref{closurespecial} of the previous subsection) and by Borwein
and Erd\'{e}lyi \cite{borweinerdelyi} and Erd\'{e}lyi
\cite{erdelyistudia} for general sequences. In all cases the
conclusion is that the elements in the closure of a nondense
M\"{u}ntz space are analytic functions. The most general result is
the following one, proved by Erd\'{e}lyi \cite{erdelyistudia} in
2003.

\begin{theorem}[Full Clarkson-Erd\H{o}s-Schwartz Theorem]\label{fullclarksonerdosschwartz}
Let $\Lambda=(\lambda_k)_{k=1}^{\infty}\subset (0,\infty)$ be a
sequence of positive real numbers such that
$M:=\Pi(\Lambda\cup\{1\})$ is a nondense M\"{u}ntz subspace of
$C[0,1]$. Then every function that belongs to the closure of $M$
in the uniform norm can be represented as an analytic function
defined on the set $\{z : z\in\mathbb{C}\setminus (-\infty,0],
\,|z|<1\}$.
\end{theorem}

\subsection{Full M\"{u}ntz theorem away from the origin}

It is remarkable that the extension of the M\"{u}ntz Theorem to
intervals away from the origin is a nontrivial task. Of course, a
linear change of variable of the form $x=bt$ allows to extend the
M\"{u}ntz Theorem to the interval $[0,b]$.
If $\Pi(\Lambda)$ is dense in $C[0,b]$, then given $f\in C[a,b]$
with $0<a<b$ one can extend $f$ with continuity to a function
$\overline{f}\in C[0,b]$ that vanishes at the origin. This
function can of course be approximated uniformly on $[0,b]$ by
elements of $\Pi(\Lambda\setminus\{0\})$, so that $f$ also belongs
to the closure of $\Pi(\Lambda\setminus\{0\})$ in $C[a,b]$. This
means that if the M\"untz condition is satisfied, then the M\"untz
polynomials are dense in $C[a,b]$.

The difficult part is to prove that the M\"{u}ntz condition is
also necessary for intervals away from the origin and, as we have
already noted, this was proved for the first time and for the
particular case of nonnegative integer exponents by Clarkson and
Erd\H{o}s. Their work was continued by L. Schwartz who proved a
full M\"{u}ntz theorem for intervals away from the origin and
general sequences of exponents. In particular, he noticed that if
we assume $0< a < b$, then the monomials $x^{\lambda }$ are
continuous functions for all $\lambda \in \Bbb{R}$ so that it
makes sense to ask for necessary and sufficient conditions for
arbitrary sequences of real numbers $\Lambda =(\lambda
_{k})_{k=0}^{\infty }\subset \mathbb{R}$ in order to make $\Pi
(\Lambda )$ a dense subspace of $C[a,b]$, and proved the following
nice result.

\begin{theorem}[Full M\"{u}ntz Theorem away from the Origin]\label{muntzawayorigin} Let
$\Lambda =(\lambda _{k})_{k=0}^{\infty }\subset \Bbb{R}$ be a
sequence of distinct real numbers, and let $0<a<b$.
Then $\Pi (\Lambda )$ is dense in $C[a,b]$ if and only if $%
\sum_{\lambda _{k}\neq 0}1/\left| \lambda _{k}\right| =\infty.$
\end{theorem}

We devote this subsection to providing a proof of this result.
Clearly, there is no loss of generality if we assume that
$0<a<b=1$. We start by assuming that the exponents can be
rearranged in such a way that they form a biinfinite sequence
$(\lambda _{k})_{k=-\infty }^{\infty }$ satisfying the following
restrictions:
\begin{itemize}
\item $\lambda _{k}>0$ for all $k>0$,

\item  $\lambda _{k}<0$ for all $k<0$,

\item $\inf_{k\in \Bbb{Z}}(\lambda _{k}-\lambda _{k-1})>0$.
\end{itemize}
We define, for each polynomial $p(z)=\sum_{|k|\leq
n}a_{k}z^{\lambda _{k}}$, the associated polynomials
\[p^{+}(z):=\sum_{0\leq k\leq n}a_{k}z^{\lambda _{k}}\qquad\text{ and}
\qquad p^{-}(z):=\sum_{-n\leq k<0}a_{k}z^{\lambda _{k}}.\]

Under these restrictions, it is possible to prove the following
relations between the uniform norms of the polynomials $p^{+}$,
$p^{-}$ and $p$:

\begin{lemma}\label{lemma+-}
Let $\Lambda = (\lambda _{k})_{k=-\infty }^{\infty }$ satisfy the
conditions we have just described and let us also assume that
$\sum_{k\in\mathbb{Z}\setminus\{0\}}1/|\lambda_k|<\infty$. Then
there exists a constant $c=c(\Lambda)$ such that
\[\|p^{+}\|_{C[a,b]}\leq
c\|p\|_{C[a,b]} \qquad\text{and}\qquad \|p^{-}\|_{C[a,b]}\leq c\|p\|_{%
C[a,b]}\]
hold for all $p\in \Pi(\Lambda)$.
\end{lemma}

\Proof.   We assume that $0<a<b=1$. It is sufficient to prove the
first inequality of the lemma since the other one is obtained from
the first via the change of variable $y=x^{-1}$. If we see the map
$p\mapsto p^{-}$ as a linear projector
${L}:\Pi (\Lambda)\rightarrow \Pi ((\lambda _{k})_{k=-\infty
}^{-1})$, the inequality we want to prove can be reformulated as:
{\sl ${L}$ is bounded whenever we use the uniform norm in the
interval
$[a,1]$ for both spaces $\Pi (\Lambda)$ and $\Pi ((\lambda
_{k})_{k=-\infty }^{-1})$}.

If ${L}$ is unbounded then there exists a sequence $%
(p_{n})_{n=0}^{\infty }\subset $ $\Pi (\Lambda)$
such that $\|p_{n}^{-}\|_{[a,1]}=1$ for all $n\geq 0$ and $%
\lim_{n\rightarrow \infty }\|p_{n}\|_{[a,1]}=0$. This clearly implies that $%
\{p_{n}^{+}\}_{n=0}^{\infty }$ is a bounded subset of $C[a,1]$
(just take into consideration that $p_{n}=p_{n}^{+}+p_{n}^{-}$ for
all $n$),
so that it is also a bounded subset of $C[0,1]$, since $%
\|p_{n}^{+}\|_{[0,1]}\leq c_{a}\|p_{n}^{+}\|_{[a,1]}$ holds for
all $n$. We can use Theorem \ref{closurespecial} and Corollary
\ref{boundedchebychevineq} to prove that there exists a
sequence of natural numbers $(n_{i})_{i=0}^{\infty }$ such that $%
(p_{n_{i}}^{+})_{i=0}^{\infty }$ converges uniformly on compact
subsets of $[0,1)$ to a certain function $f^{+}=$
$\sum_{k=0}^{\infty }a_{k}z^{\lambda_{k}}$ analytic on
$\mathbb{D}_{1}:=(\Bbb{C}\setminus (-\infty ,0])\cap
\mathbb{D}(0,1)$, and the sequence
$(p_{n_{i}}^{-})_{i=0}^{\infty}$ converges uniformly on compact
subsets of $(a,\infty )$ to a certain function $f^{-}=$
$\sum_{k=-\infty }^{-1}a_{k}z^{\lambda _{k}}$ which is analytic on
$\Bbb{E}_{a}:=(\Bbb{C}\setminus (-\infty ,0])\cap
(\Bbb{C}\setminus \overline{\mathbb{D}}(0,a))$. The last claim can
be proved by just making the change of variable $t=a/x$, since
then $r_{n}^{+}(t)=p_{n}^{-}(a/x)\in \Pi
((-\lambda _{k})_{k=-\infty }^{-1})$ is a Cauchy sequence in $C%
[a,1]$, so that we can assume that
$(r_{n_{i}}^{+})_{i=0}^{\infty }$ converges uniformly on compact
subsets of $[0,1)$ to a certain $h^{+}=$ $\sum_{k=-\infty
}^{-1}h_{k}z^{-\lambda _{k}}$ which is analytic in
$\mathbb{D}_{1}:=(\Bbb{C}\setminus (-\infty ,0])\cap
\mathbb{D}(0,1)$, for an
adequate choice of $(n_{k})_{k=0}^{\infty }$. We then go back with
the change of variable, obtaining that
$f^{-}=h^{+}(a/z)=\sum_{k=-\infty }^{-1}a_{k}z^{\lambda _{k}}$
(where $a_{k}=h_{k}a^{-\lambda _{k}}$ for all $k$) is analytic in
$\Bbb{E}_{a}$ and
$(p_{n_{k}}^{-})_{k=0}^{\infty }$
converges uniformly to $f^{-}$ on compact subsets of $(a,\infty )$. Now, $%
\lim_{n\rightarrow \infty }\|p_{n_{k}}\|_{[a,1]}=0$ and $%
p_{n_{k}}=p_{n_{k}}^{+}+p_{n_{k}}^{-}$ for all $i$, so that
$f^{+}+f^{-}=0$ in $(a,1)$. This implies that
\[
g(z):=\left\{
\begin{array}{ccc}
f^{+}(\ee^{z}), &  & {\rm Re}(z)<0, \\[5pt]
-f^{-}(\ee^{z}), &  & {\rm Re}(z)>\log a,
\end{array}
\right.
\]
can be extended as a bounded entire function \cite[page
181]{borweinerdelyi}. It follows from Liouville's theorem
that $g={\rm const}$, so that $g=0$ since $\lim_{t\rightarrow \infty }f^{-}(t)=0$. Hence $%
f^{+}=0$ in $[0,1)$ and $f^{-}=0$ in $(a,\infty )$, which implies $%
\lim_{k\rightarrow \infty }\|p_{n_{k}}^{-}\|_{[a,1]}=0$, a
contradiction.\endproof

We now characterize the closure of $\Pi (\Lambda)$ whenever
$\Lambda = (\lambda _{k})_{k=-\infty }^{\infty }$ satisfies the
additional condition $\sum_{\lambda _{k}\neq 0}1/\left| \lambda
_{k}\right| <\infty $.

\begin{theorem}\label{teo22}
Assume that $\Lambda = (\lambda _{k})_{k=-\infty }^{\infty }$ satisfies $%
\lambda _{k}>0$ for all $k>0$, $\lambda _{k}<0$ for all $k<0$,
$\inf_{k\in
\Bbb{Z}}(\lambda _{k}-\lambda _{k-1})>0$, and \ $\sum_{%
\lambda _{k}\neq 0}1/\left| \lambda _{k}\right|<\infty $. Then the
elements of the closure of $\Pi (\Lambda)$ in $C[a,b]$ can be
extended as analytic functions to the domain
\[
\{z : z\in \Bbb{C}\setminus (-\infty ,0],\, a<|z|<b\}\text{.}
\]
\end{theorem}

\Proof. Let us assume, without loss of generality, that $0<a<b=1$.
We have already proved in Lemma \ref{lemma+-}, under the
hypotheses we have imposed on
$\Lambda$, that if $f$ belongs to the closure of $\Pi
(\Lambda)$ in $C[a,b]$ then $%
f=f^{+}+f^{-}$, where $f^{+}$ is analytic in $\mathbb{D}_{1}$ and
$f^{-}$ is analytic in $\Bbb{E}_{a}$. Hence $f$ is analytic in
\[
\mathbb{D}_{1}\cap \Bbb{E}_{a}=\{z : z\in \Bbb{C}\setminus (-\infty ,0],\, a<|z|<1\}%
\text{,}
\]
and the proof is complete.\endproof

\Proof {the Full M\"{u}ntz theorem away from the Origin}. Let us
decompose the proof into the following four cases:

\textbf{Case 1. }$\{\lambda _{k}\}_{k=0}^{\infty }$ has some
accumulation point $\lambda \neq 0$.

We can assume without loss of generality that $\lambda>0$.
(Otherwise, consider the map $S:C[a,b]\to C[a,b]$ given by
$S(f)(x)=x^{-\lambda+1}f(x)$ and take into account that $S$ is a
linear isomorphism of Banach spaces, so that it transforms dense
subspaces into dense subspaces and {\it vice versa}.) Hence this
case is an easy corollary of the Full M\"{u}ntz Theorem for the
interval $[0,b]$.

\textbf{Case 2. }$0$ is an accumulation point of $\{\lambda
_{k}\}_{k=0}^{\infty }$.

This case is reduced to Case 1 as follows: first we consider the
sequence $(\lambda _{k}+1)_{k=0}^{\infty }$, which is in Case
1, so that $\Pi ((\lambda _{k}+1)_{k=0}^{\infty })$ is dense in
$C[a,b]$. Now we consider the isomorphism of Banach spaces given
by
$T:C%
[a,b]\rightarrow C[a,b]$, $(Tf)(x):=x^{-1}f(x)$ and
conclude that $\Pi ((\lambda _{k})_{k=0}^{\infty })$ is dense in
$C[a,b]$, since $T(\Pi ((\lambda _{k}+1)_{k=0}^{\infty }))=\Pi
((\lambda _{k})_{k=0}^{\infty})$.

\textbf{Case 3. }$\{\lambda _{k}\}_{k=0}^{\infty }$ has no
accumulation
points and $\sum_{\lambda _{k}>0}\frac{1}{\lambda _{k}}%
=\infty $ or $\sum_{\lambda _{k}<0}\frac{1}{\left| \lambda
_{k}\right| }=\infty $.

In this case, we may assume, without loss of generality, that
$0\notin $ $\{\lambda _{k}\}_{k=0}^{\infty }$ since otherwise we
can take $\varepsilon
>0$ such that $0\notin \{\lambda _{k}+\varepsilon
\}_{k=0}^{\infty }$ and (by using the same kind of arguments as in
Case 2) we have that the density of $\Pi ((\lambda
_{k})_{k=0}^{\infty })$ and the density of $\Pi ((\lambda
_{k}+\varepsilon )_{k=0}^{\infty })$ are equivalent claims.

If $\sum_{\lambda _{k}>0}\frac{1}{\lambda _{k}}=\infty $, then we
can use the Full M\"{u}ntz Theorem on $C[0,b]$ to
conclude the proof. If $\sum_{\lambda _{k}<0}\frac{1}{%
|\lambda _{k}|}=\infty $ but $\sum_{\lambda _{k}>0}\frac{1}{%
\lambda _{k}}<\infty $, then we use the change of variable
$t=1/x$, and that $S:C[a,b]\rightarrow C[1/b,1/a]$,
$S(f)(x):=f(1/x)$ is a linear isometry
($\|S(f)\|_{C[1/b,1/a]}=\|f\|_{C[a,b]}$ is clear), to
prove that $\Pi ((\lambda _{k})_{k=0}^{\infty })$ is dense in $C%
[a,b]$ if and only if $S(\Pi ((\lambda _{k})_{k=0}^{\infty }))=$
$\Pi ((-\lambda _{k})_{k=0}^{\infty })$ is dense in
$C[1/b,1/a]$, which puts us once more in the case $\sum_{\lambda
_{k}>0}\frac{1}{\lambda _{k}}=\infty $.

\textbf{Case 4. }$\sum_{\lambda _{k}\neq 0}1/\left| \lambda
_{k}\right| <\infty $.

We rearrange the sequence $(\lambda _{k})_{k=0}^{\infty }$ as
$(\lambda _{k}^{\ast })_{k=-\infty }^{\infty }=(\lambda
_{k})_{k=0}^{\infty }$, with $\lambda _{k}^{\ast }<\lambda
_{k+1}^{\ast }$ for all $k\in \Bbb{Z}$, $\lambda _{k}^{\ast }<0$
if $k<0$ and $\lambda _{k}^{\ast }>0$ if $k>0$. Then there exists
a sequence $\Gamma:=(\gamma_{k})_{k=-\infty }^{\infty }$ such that
$\inf_{k\in\Bbb{Z}}(\gamma_{k}-\gamma _{k-1})>0$, $\sum_{k\in
\Bbb{Z}}1/|\gamma _{k}|<\infty $, and $\gamma _{k}<\gamma _{k+1}$,
$|\gamma _{k}|<|\lambda _{k}^{\ast }|$ for all $k\in \Bbb{Z}$,
$\gamma _{k}<0$ if $k<0$ and $\gamma _{k}>0$ if $k\geq 0$. Now, it
follows from Theorem \ref{teo22} that there is an $m$ such that
$x^m\notin \overline{\Pi (\Gamma)}$
and from the comparison theorem for real exponents (see also
Remark \ref{remarkoncomparisontheorem}) that $x^m\notin
\overline{\Pi
(\Lambda)}$.

This completes the proof of the Full M\"{u}ntz Theorem away from
the origin.\endproof

\subsection{Full M\"{u}ntz theorem for measurable sets}

Borwein and Erd\'elyi have recently published several papers in
which they prove a Full M\"{u}ntz Theorem for the spaces $C(A)$
and
\[L_{w}^{q}(A)=\{f:\left(\int_{A}|f(x)|^{p}w(x)\dd x\right)^{1/p}<\infty\},\]
for sets $A$ with positive Lebesgue measure and weight functions
$w$ (i.e., $w>0$ is measurable in the sense of Lebesgue). To be
more precise, we state here one of their main results (see
\cite{borwein6},\cite{borwein5}, and \cite{borwein4}).

\begin{theorem}[Borwein-Erd\'elyi]\label{muntzmeasurable}
If $\Lambda =(\lambda _{k})_{k=-\infty }^{\infty
}\subset \Bbb{R}$ is a sequence of distinct real numbers with $%
\lambda _{k}<0$ for all $k<0$, $\lambda _{k}\geq 0$ for all $k\geq
0$ such that $\sum_{\lambda _{k}\neq 0}1/|\lambda_{k}|<\infty$ and
$A\subset (0,\infty )$ is a set with positive Lebesgue measure
such that $\inf A>0$, then $\Pi (\Lambda )$ is not dense in
$L_{w}^{q}(A)$ for all weight functions $w:A\rightarrow \lbrack
0,\infty )$ with $\int_{A}w>0$ and all $q\in (0,\infty )$.
Moreover, every function that belongs to the closure of \ $\Pi
(\Lambda )$ in $L_{w}^{q}(A)$ can be analytically extended to the
domain $\{z: z\in \Bbb{C}\setminus (-\infty ,0],\, a_{w}<\left| z\right| <b_{w}\}$%
, where
\begin{eqnarray*}
a_{w} &:=&\inf \{y\in \lbrack 0,\infty ):\int_{A\cap (0,y)}w>0\}, \\[10pt]
b_{w} &:=&\sup \{y\in \lbrack 0,\infty ):\int_{A\cap (y,\infty
)}w>0\}.
\end{eqnarray*}
Finally, if
\[
\inf \{\lambda _{k}-\lambda _{k-1}:k\in \Bbb{Z}\}>0
\]
then all functions that belong to the closure of $\Pi (\Lambda )$ in $%
L_{w}^{q}(A)$ admit a representation of the form
\[
f(x)=\sum_{k=-\infty }^{\infty }a_{k}x^{\lambda _{k}}, \qquad x\in
A\cap (a_{w},b_{w}).
\]
\end{theorem}

The key idea  in proving Theorem \ref{muntzmeasurable} is to use
Egorov's theorem and the following important polynomial inequality
(see \cite {borwein5}):

\begin{theorem}[Remez Inequality for M\"{u}ntz Polynomials]
Let $\Lambda =(\lambda _{k})_{k=-\infty }^{\infty }\subset
\Bbb{R}$ be an arbitrary sequence of real numbers. If\/
$\sum_{\lambda_k\neq 0}1/\left| \lambda _{k}\right|<\infty$, then
for all sets $A\subset \lbrack 0,\infty )$ with Lebesgue measure
$m(A)>0$ and all intervals $[\alpha ,\beta ]\subset ({\rm
ess\,inf} (A), {\rm ess\,sup} (A))$, there exists a constant
$c=c(\Lambda ,A,\alpha ,\beta )$ such that
\[
\left\| p\right\| _{C[\alpha ,\beta ]}\leq c\left\| p\right\| _{C%
(A)}
\]
for all $p\in \Pi (\Lambda )$.
\end{theorem}

In fact, it follows from this theorem that we can easily prove the
following result, which is a main step in the proof of the
corresponding Full M\"{u}ntz Theorem for sets of positive Lebesgue
measure:

\begin{corollary}
Let $\Lambda =(\lambda _{k})_{k=-\infty }^{\infty }\subset
\Bbb{R}$ be an arbitrary sequence of real numbers such that $
\sum_{\lambda_k\neq 0}1/\left| \lambda _{k}\right| <\infty$. Then,
for any set $A\subset \lbrack 0,\infty )$ with positive Lebesgue
measure $m(A)>0$, we have that if the sequence of
polynomials $(p_{n})_{n=0}^{\infty}\subset \Pi (\Lambda )$
converges pointwise to $f\in C(A)$, then for all $[\alpha ,\beta
]\subset (a,b):=({\rm ess\,inf}(A),{\rm ess\,sup}(A))$,
$(p_{n})_{n=0}^{\infty}$ is a Cauchy sequence in $C[\alpha ,\beta
]$.
\end{corollary}

\Proof. First of all, we would like to recall
that Egorov's theorem guarantees that if $(f_{n})$ is a sequence
of measurable functions on $A$ (where $0<m(A)<\infty $) that
converges almost everywhere to a certain function $f$ (that is
finite almost everywhere on $A$), then for all $\varepsilon >0$
there exists a measurable set $B\subset A$ such that $m(A\setminus
B)<\varepsilon $ and $(f_{n})$ converges uniformly on $B$ to $f$.

Let $(p_{n})_{n=0}^{\infty}$ and $f$  satisfy the hypotheses of
this corollary and let $[\alpha ,\beta ]\subset (a,b)$. It follows
from the definition of $(a,b)$ and from Egorov's theorem that
there are sets of positive Lebesgue measure
\[
B_{1}\subset A\cap (0,\alpha )\text{ and }B_{2}\subset A\cap
(\beta ,\infty )
\]
such that $(p_{n})_{n=0}^{\infty}$ converges uniformly to $f$ on
$B=B_{1}\cup B_{2}$.

Now, the application of the Remez inequality for Müntz polynomials
on $[\alpha ,\beta ]\subset (\rho ,\sigma )$ (where $\rho :={\rm
ess\,inf} (B)$ and $\sigma :={\rm ess\,sup} (B)$),
\[
\left\| p_{i}-p_{j}\right\| _{[\alpha ,\beta ]}\leq C(B,[\alpha
,\beta ],\Lambda )\left\| p_{i}-p_{j}\right\| _{B}\text{,}
\]
proves that $(p_{n})_{n=0}^{\infty}$ is a Cauchy sequence in
$C[\alpha ,\beta ]$.\endproof

\subsection{Full M\"{u}ntz theorem for countable compact sets}
It is quite surprising that for a long period of time the
M\"{u}ntz Theorem has been studied in many cases but not for the
space $C(K)$ with  $K$ a countable compact set. This is surprising
because, in principle, this case should be the easiest one. This
question has been addressed quite recently by the author
\cite{almirabms}, and it transpires that in many cases the
M\"{u}ntz condition can be weakened in a sensible way when dealing
with countable compact sets. In particular, the following result
holds.

\begin{theorem}[Almira, 2006]\label{teoalmira} Let $K\subset [0,\infty)$ be an
infinite countable compact set and let
$\Lambda=(\lambda_k)_{k=0}^{\infty}\subset \mathbb{R}$ be a fixed
sequence of exponents, satisfying $\lambda_0=0$. Then the
following holds:
\begin{description}
\item[i)] If $\Lambda\subset [0,\infty)$ is an infinite bounded
sequence and $K\setminus \{0\}$ is compact then $\Pi (\Lambda)$ is
dense in $C(K)$. \item[ii)] If $\Lambda \subset [0,\infty)$ and
$K$ does not contain strictly increasing infinite sequences then
$\Pi (\Lambda)$ is dense in $C(K)$ if and only if
$\#\Lambda=\infty$. Moreover, if $\Lambda \subset (-\infty,0]$ and
$K$ does not contain strictly decreasing  infinite sequences then
$\Pi (\Lambda)$ is dense in $C(K)$ if and only if
$\#\Lambda=\infty$.
\end{description}
\end{theorem}

\Proof.  The main idea in the proof is to use the Riesz
Representation Theorem. Clearly, the unique measures that exist
for countable compact sets are atomic. Thus, if
$K=\{0\}\cup\{t_i\}_{i=1}^{\infty}$ then $L\in C^*(K)$ if and only
if $L(f)=\alpha_0f(0)+\sum_{i=1}^{\infty}\alpha_if(t_i)$ for a
certain sequence $(\alpha_i)_{i=0}^{\infty}$ such that
$\sum_{i=0}^{\infty}|\alpha_i|<\infty$. Thus, as a consequence of
the Hahn-Banach Theorem, ${\rm
span\,}\{x^{\lambda_k}\}_{k=0}^{\infty}$ is dense in $C(K)$ if and
only if the following holds: if $\sum_{i=0}^{\infty}\alpha_i=0$,
\[
\sum_{i=1}^{\infty}\alpha_it_i^{\lambda_k}=0, \quad k=1,2,\ldots,
\quad \text{ and } \quad \sum_{i=0}^{\infty}|\alpha_i|<\infty,
\]
then $\alpha_i=0$ for all $i\geq 0$.

Thus, let us assume that
\[\sum_{i=0}^{\infty}\alpha_i=0,\quad
\sum_{i=1}^{\infty}\alpha_it_i^{\lambda_k}=0, \quad
k=1,2,\ldots,\quad \mbox{\rm and}\quad
\sum_{i=0}^{\infty}|\alpha_i|<\infty.\]
 Then we set
$\Gamma:=\{t_i:i\geq 1, \alpha_i\neq 0\}$
 and we take $\gamma:=\sup\Gamma$.
Clearly, $\gamma\in K$ since $K$ is compact. If $\Gamma=\emptyset$
then $L(f)=\alpha_0f(0)$ and $L(1)=0$ implies $\alpha_0=0$, which
ends the proof. If $\Gamma\neq \emptyset$ then $\gamma>0$ and
there exists $t_s\in K$ such that $\gamma=t_s$. Thus, we take
$t_a\in K$ such that $t_a<t_s$ and we set
$z_{\lambda}:=(t_a/t_s)^{\lambda}$. Clearly, the equation
$z_{\lambda}^{p_j}=(t_j/t_s)^{\lambda}$ is uniquely solved by
$p_j=(\ln(t_j/t_s))\bigr/ \ln (t_a/t_s)$, which is a positive real
number for all $j\neq s$ . Hence $L(x^{\lambda_k})=0$,
$k=0,1,2,\ldots$, can be written in the following equivalent way:
\[
0=\sum_{i=0}^{\infty}\alpha_i\quad \text{ and } \quad 0
=(t_s)^{\lambda_k}\sum_{t_i\in\Gamma}\alpha_i\left(\frac{t_i}{t_s}\right)^{\lambda_k},
\quad k=1,2,\ldots .
\]
Hence $\varphi(z_{\lambda_k})=0$ for all $k\geq 1$, where
\[
\varphi(z):=\sum_{t_i\in\Gamma}\alpha_iz^{p_i} .
\]

We decompose the proof into several steps, according to the
boundedness properties of the sequence of exponents $\Lambda$.

\noindent \textbf{Step 1.} $\Lambda\subset [0,\infty)$ and
$\lim_{k\to\infty}\lambda_k=\infty$ and $K$ does not contain
strictly increasing infinite sequences.

Under these conditions, it is clear that $t_s\in \Gamma$ and
$\lim_{k\rightarrow \infty }z_{\lambda_k}=0$. Thus
$\varphi(0)=\lim_{k\rightarrow\infty}\varphi(z_{\lambda_k})\hfill\newline=0$
since $\varphi(z)$ is continuous at the origin. On the other hand,
$t_s\in \Gamma$ implies that $\alpha_s\neq 0$. Hence we can use
the fact that $\varphi(z)=
\sum_{t_i\in\Gamma\setminus\{t_s\}}\alpha_iz^{p_i}+\alpha_s$
(since $p_s=(\ln 1)/\ln(t_a/t_s)=0$) to claim that
$\varphi(0)=\alpha_s\neq 0$, a contradiction.

\noindent \textbf{Step 2.}
$\Lambda=(\lambda_k)_{k=0}^{\infty}\subset [0,\infty)$ is bounded,
$0\neq \lim_{k\to\infty}\lambda_k$.

Clearly, we can assume without loss of generality that $\Lambda$
is  itself a convergent sequence. We note that
$\varphi(z)=\sum_{t_i\in\Gamma}\alpha_iz^{p_i}$ is analytic in the
open set $\Omega=\{z:|z|<1,|1-z|<1\}$. If
$\lim_{k\to\infty}\lambda_k=\lambda^*\neq 0$ then
$\lim_{k\to\infty}z_{\lambda_k}=z_{\lambda^*}\in (0,1)\subset
\Omega$. Hence $\varphi$ vanishes on a set with accumulation
points inside $\Omega$, so that $\varphi(z)$ vanishes identically
on $\Omega$ and $\alpha_i=0$ for all $i>0$. If $0\not \in K$ the
proof is complete. On the other hand, if $0\in K$ then
$0=L(1)=\sum_{t_i\in \Gamma}\alpha_i+\alpha_0=\alpha_0$ and the
proof is also complete.

\noindent \textbf{Step 3.} $\lim_{k\to\infty}\lambda_k=0$ and
$K\setminus\{0\}$ is compact.

In this case, we can use the following trick: the equations
\[
0=\sum_{t_i\in\Gamma}\alpha_it_i^{\lambda_k}, \quad k=1,\ldots
\]
can be rewritten as
\[
0=\sum_{t_i\in\Gamma}\beta_it_i^{\lambda_k^*}, \quad k=1,\ldots,
\]
where $\beta_i:=\alpha_i/t_i$ for all $i$ and
$\lambda_k^*:=\lambda_k+1$ for all $k$ (taking into account that
$\sum_{t_i\in\Gamma}|\beta_i|<\infty$ since $K\setminus \{0\}$ is
compact). Thus $\lim_{k\to\infty}\lambda_k^*=1$ and we conclude
that $\alpha_j/t_j=0$ for all $j$. The proof follows.

\noindent \textbf{Step 4.} $\Lambda\subset \mathbb{R}$ and
$K\setminus\{0\}$ is compact.

Clearly, if $\Lambda$ is an infinite set then it contains either
infinitely many positive elements or infinitely many negative
elements. Thus, we may assume that either $\Lambda\subset
[0,\infty)$ or $\Lambda\subset (-\infty,0]$. The first case has
been already studied in Steps 1 and 2. Thus, let us assume that
$\Lambda\subset (-\infty,0]$ and
$L(f)=\alpha_0f(0)+\sum_{j=1}^{\infty}\alpha_jf(t_j)\in C^*(K)$.
Then the equations $L(x^{\lambda_k})=0$, $k=0,1,\ldots$, can be
rewritten as
\[
\sum_{i=0}^{\infty}\alpha_i=0 \qquad \text{ and }\qquad
\sum_{i=1}^{\infty}\alpha_i(1/t_i)^{\lambda_k}=0, \quad
k=1,2,\ldots .
\]
This means that the functional given by
\[
S(f):=\alpha_0f(0)+\sum_{j=1}^{\infty}\alpha_jf(1/t_j),
\]
which belongs to $C^*(E)$, where
$E:=\{0\}\cup\{1/t_j\}_{j=1}^{\infty}$, which is a countable
compact subset of $[0,\infty)$ since $K\setminus\{0\}$ is compact,
satisfies $S(x^{-\lambda_k})=0$ for all $k\geq 0$. Moreover, if
$K$ does not contain decreasing sequences then $E$ does not
contain increasing sequences. Now, we use the results proved in
Steps 1, 2 and 3 to conclude that $\alpha_i=0$ for all $i$.
\endproof

\begin{remark} \label{remark55}
{\rm \ There is another proof of step 2. Taking into consideration that
$t_i^{\lambda_k}=\exp(\lambda_k\log t_i)$ for all $i$, we have
that the relations
\[
\sum_{t_i\in \Gamma}\alpha_it_i^{\lambda_k}=0, \qquad k=1,2,\ldots
\]
are equivalent to the relations
\[
\Psi(\lambda_k)=0, \qquad k=1,2,\ldots
\]
where
\[
\Psi(z):=\sum_{t_i\in\Gamma}\alpha_i\exp((\log t_i)z)
\]
is an entire function of exponential type. This means, in
particular, that $\Lambda$ cannot be an infinite bounded sequence
(otherwise, $\Psi$ should vanish everywhere).
}
\end{remark}

\begin{remark} \label{remark56} {\rm \ Clearly, if $\Lambda$ is bounded then
$|\lambda_k|^{-1}\geq 1/\sup \Lambda $ for all $\lambda_k\neq 0$.
Hence $\sum_{k=1}^{\infty}|\lambda_k|^{-1}=\infty$ and case (i) of
Theorem \ref{teoalmira}  follows from the M\"{u}ntz Theorem away
from the origin $($Theorem \ref{muntzawayorigin}$)$ whenever
$0\not \in K$. This proof uses a very difficult result in order to
prove a simpler one. This is the reason we gave our own elementary
proof of this fact.
}
\end{remark}

\begin{remark} \label{remark57} {\rm \ There are many countable compact sets
with the property that they do not have (strictly) increasing
sequences. An interesting example is given by:
\[
K=\{0\}\cup \{1/n\}_{n=1}^{\infty}
\cup\{1/n+1/m\}_{n,m=1}^{\infty}\;.
\]
Obviously, this compact set has infinitely many accumulation
points and it has no increasing sequences! These cases are covered
by Theorem \ref{teoalmira} above.
}
\end{remark}

\noindent \textbf{Open question.}  We have already shown that in
order to give a Full M\"{u}ntz Theorem for the general case (i.e.,
for arbitrary countable compact sets $K\subset [0,\infty)$), it is
a good idea to study the zero sets of the M\"{u}ntz type series
\[
\varphi(z)=\sum_{i=1}^{\infty}\alpha_jz^{p_j},
\]
 where $(p_j)_{j=1}^{\infty}$ decreases to zero,
$\sum_{j=1}^{\infty}|\alpha_j|<\infty$ and, for the case in which
$K\setminus\{0\}=\{t_i\}_{i=1}^{\infty}$ is compact, the zero sets
of the entire functions of exponential type given by
\[
\Psi(z)=\sum_{t_i\in\Gamma}\alpha_i\exp((\log t_i)z) \qquad \text{
where } \quad \sum_{j=1}^{\infty}|\alpha_j|<\infty.
\]
Is it possible to find a series $\varphi(z)$ with a sequence of
infinitely many zeros $(z_k)_{k=0}^{\infty}$ that converges to
zero? What about a function $\Psi(z)$ with infinitely many zeros?
These questions seem to be still open and not easy to solve.

\medskip
\noindent \textbf{Acknowledgement. } The author is infinitely
grateful to the referee of a previous version of this paper. With
his (her) help the manuscript improved not only its readability
but also the details in many proofs. Thanks.

\def\atextbf#1{#1}

\bigskip

J. M.\ Almira.

Departamento de Matem\'{a}ticas. Universidad de Ja\'{e}n.

E.U.P. Linares 23700 Linares (Ja\'{e}n) Spain

email: jmalmira@ujaen.es

\endddoc